\documentclass[11pt,leqno]{article}
\usepackage{amsmath,amsthm,amssymb}
\usepackage{amscd}
\usepackage[matrix,arrow]{xy}

\author{Florin Ambro}

\newcommand{\N}{{\mathbb N}}
\newcommand{\Q}{{\mathbb Q}}
\newcommand{\R}{{\mathbb R}}
\newcommand{\Z}{{\mathbb Z}}

\newcommand{\cG}{{\mathcal G}}
\newcommand{\cH}{{\mathcal H}}
\newcommand{\cI}{{\mathcal I}}

\newcommand{\cO}{{\mathcal O}}

\newcommand{\cN}{{\mathcal N}}
\newcommand{\cS}{{\mathcal S}}

\newcommand{\Supp}{\operatorname{Supp}}

\newcommand{\mult}{\operatorname{mult} }
\newcommand{\codim}{\operatorname{codim}}
\newcommand{\Pic}{\operatorname{Pic}}
\newcommand{\im}{\operatorname{Im}}
\newcommand{\Coker}{\operatorname{Coker}}

\newtheorem{thm}{Theorem}[section]  
\newtheorem{lem}[thm]{Lemma}  
\newtheorem{prop}[thm]{Proposition}   
\newtheorem{cor}[thm]{Corollary}  
\newtheorem{defn}{Definition}[section] 
\newtheorem{rem}{Remark}[section] 
\newtheorem{example}{Example}[section] 
 
\newtheorem{conj}{Conjecture}

\begin{document}

\title{The Adjunction Conjecture
                     and its applications}
\date{}
\maketitle

\pagestyle{plain}
\pagenumbering{roman}
\setcounter{page}{2}

\begin{center}

 {\bf Abstract}

\end{center}

Adjunction formulas are fundamental tools in the
classification theory of algebraic varieties.
In this paper we discuss adjunction formulas for fiber spaces 
and embeddings, extending the known results
along the lines of the Adjunction Conjecture, independently 
proposed by Y. Kawamata and V. V. Shokurov.
\par
 As an application, we simplify Koll\'ar's proof for the Anghern 
and Siu's quadratic bound in the Fujita's Conjecture. We also 
connect adjunction and its precise inverse to the problem of building 
isolated log canonical singularities.

\tableofcontents

\pagenumbering{arabic}
\pagestyle{plain}

\section*{Introduction}
 The rough classification of projective algebraic varieties attempts to 
divide them according to properties of their canonical class $K_X$. 
Therefore, whenever two varieties are closely 
related, it is essential to find formulas comparing their canonical 
divisors. Such formulas are called {\em adjunction formulas}. We first
look at the following examples:

\begin{itemize}
  \item[(1)] Let $C$ be a smooth curve on a smooth projective
surface $X$. Then $$(K_X+C)|_C \sim K_C $$ 

  \item[(2)] Let $C \subset {\mathbb P}^2$  be the curve defined by the 
equation 
$x^2 z-y^3=0$, with normalization $C^{\nu} \simeq {\mathbb P}^1$. 
Then $$(K_{{\mathbb P}^2}+C)|_{C^\nu} \sim
K_{C^{\nu}}+2\cdot P$$ where 
$P$ is the point of 
$C^\nu$ above the cusp $(0:0:1) \in C$. 
  \item[(3)] Let $X$ be a nonsingular variety and $W \subset X$ be 
the complete intersection of the nonsigular divisors $D_1, \ldots, D_k$.
Then $$(K_X+\sum_i D_i)|_W \sim K_W $$ 
\end{itemize}

The appearance of the divisor $B_{C^\nu}=2\cdot P$, 
called the {\em different} of 
the log divisor $K_X+C$ on $C^\nu$, is what makes the adjunction formula 
symmetric (the different happens to be trivial in the other two examples).
This phenomena is called {\em subadjunction} in \cite[5-1-9]{KMM} and 
it was first observed by Miles Reid in this context. 
\par
The examples above suggest that a {\em log divisor} $K_X+B$ has a 
natural {\em residue} $K_W+B_W$ on the normalization $W$ of the
intersection of components of $B$ with coefficient $1$. 
It is also expected that the {\em moduli part} \cite{Ka2,Ka3} 
$M_W=(K_X+B)|_W-(K_W+B_W)$ is {\em semiample}, that is one of its
multiples is base point free. 
These ideas are formalized in the following conjecture, independently 
proposed by Yujiro Kawamata and Vyachelsav Shokurov.

\begin{conj}[The Adjunction Conjecture] Let $(X,B)$ be a log variety and
let $j:W \to X$ be the normalization of an irreducible component
of $LCS(X,B)$ such that $K_X+B$ is log canonical in the generic point
of $j(W)$. There exists a canonically defined $\R$-divisor 
$B_W$ on $W$, called the {\em different} of the log divisor $K+B$ on 
$W$, with the following properties:
\begin{enumerate}
   \item $(W,B_W)$ is a log variety, that is $K_W+B_W$ is $\R$-Cartier
         and $B_W$ is an effective divisor;
   \item The induced map $j:(W,B_W) \to (X,B)$ is {\em log proper},
         i.e. for each closed subvariety $Z \subset W$, there 
         exists a natural number $N \in \N$ such that
         $$\frac{a(j(Z);X,B)}{N}  \le a(Z;W,B_W) \le a(j(Z);X,B)$$
   \item (Freeness) $(K+B)|_W \sim_{\R} K_W+B_W+M$, where $M$ is an 
         $\R$-free divisor on $W$;
   \item (Boundness) If $K+B$ is $\Q$-Cartier of index $r$, 
         there exists a natural number $$b=b(r,dimW,dimX) \in \N$$ such 
         that $b(K_W+B_W)$ and $b M$ are Cartier divisors.
\end{enumerate}
\end{conj}

 Chapter $1$ is introductory. We introduce in Chapter $2$
Shokurov's {\em minimal log discrepancies} $a(Z;X,B)$, measuring
the singularities of a log pair $(X,B)$ in a closed subvariety $Z$,
and we discuss two conjectures: the lower semicontinuity of
minimal log discrepancies and the precise inverse of 
adjunction (the exceptional case). 
\par
Chapters $3$ and $4$ are an expanded 
version of Y. Kawamata's papers on adjunction \cite{Ka1,Ka2,Ka3}.
In Chapter $3$ we define the {\em discriminant} $B_Y$ of a log
divisor $K_X+B$ along a morphism $f:X \to Y$. It measures
the singularities of $K+B$ above the codimension $1$ points 
of $Y$. The discriminant appears in \cite{Ka2,Ka3} for special morphisms, 
as well as in \cite[5.12,9.12]{Mo}, where it is called the 
{\em negligible part}. For instance, a fiber space of smooth varieties 
$f:X \to Y$ with simple normal crossing ramification divisor is semistable 
in codimension $1$ iff the discriminant of $K_X$ is trivial. 
\newline
We prove the finite base change formula for discriminants 
and we propose the Base Change Conjecture, claiming that 
the {\em birational} base change formula for discriminants holds 
for {\em log Calabi-Yau fiber spaces}. 
The Base Change Conjecture is intuitively equivalent to the 
log properness property in the Adjunction Conjecture.
 Finally, we present an extension of the result of Y. Kawamata on
the nefness of the moduli part for certain log Calabi-Yau fiber spaces.
The Base Change Conjecture implies that Kawamata's positivity result 
holds for every log Calabi-Yau fiber space. 
\par
 In Chapter $4$ we introduce the {\em different} of a log
divisor $K_X+B$ on a log canonical (lc) center $W \subset X$ such that 
$a(\eta_W;B)=0$. We restrict to the case when $W$ is an {\em exceptional} 
lc center, the higher codimensional equivalent of generic
pure log terminality. For instance, any codimension $1$ lc center
is exceptional. However, the definition and the properties of
the different hold for non-exceptional lc centers too, once the 
basic adjunction calculus is extended from normal to seminormal varieties.
 The different is the discriminant of a log Calabi-Yau fiber space, 
so the properties of the latter translate into properties of the 
former. We also show that the Base Change Conjecture reduces
the first two properties of the different stated in the 
Adjunction Conjecture to the case $codim(W,X)=1$.
\par
 We conclude this chapter with an extension of Kawamata's adjunction 
formula \cite{Ka3}. This weak version of adjunction is enough for
certain applications. We use it in Chapter $5$ to reobtain the known 
quadratic bound for building isolated log canonical singularites, found 
by U. Anghern and Y.T. Siu in the analytic case, later adapted by 
J. Koll\'ar to the algebraic case.
\par
Finally, we discuss applications of adjunction,
as an excellent tool for inductive arguments in Higher Dimensional 
Algebraic Geometry.
Chapter $5$ deals with the problem of {\em building log canonical 
singularities}. If $x \in X$ is a closed point such that $a(x;B) \ge 0$
we search for effective $\Q$-divisors $D$ such that 
$a(x;B+D)=0$. To make the 
problem nontrivial, we fix an ample $\Q$-Cartier divisor $H$
on $X$ and ask what is the infimum $bld_x(B;H)$ of all $c>0$ 
for which there exists a divisor $D \sim_\Q cH$ with the above property.
\par
 If $X$ is a curve, then $bld_x(B;H)=a(x;B)/deg_X(H)$, in other words
$bld_x(B;H)$ is controlled by the (global) numerical properties of
$H$ and the (local) invariants of the singularity of the log
variety $(X,B)$ at $x$. 
The optimal bound for $bld_x(B;H)$ is stated as Conjecture~\ref{bld_conj}. 
We show that the first two properties stated in the Adjunction Conjecture 
reduce Conjecture~\ref{bld_conj} to its inductive step, stated as 
Conjecture~\ref{bld_step}. We also provide some evidence for 
Conjecture~\ref{bld_step}. 
\newline 
A lemma of Y. Kawamata translates any upper bound
for $bld_x(B;H)$ in effective results on the global generation 
of (log) adjoint line bundles on projective varieties. In particular,
the first two properties stated in the Adjunction Conjecture 
and Conjecture~\ref{bld_step} imply a stronger version of
Fujita's Conjecture.
\par
 As a final remark, Chapter $4$ leaves the Adjunction Conjecture still 
in a hypothetical form. The only satisfactory case so far 
is $codim(j(W),X)=1$, where all the properties are checked, with the 
exception of the {\em precise inverse of adjunction}. 
However, assuming that $j(W)$ is exceptional and that 
the Base Change Conjecture holds true we can summarize the 
known results as follows:
\begin{itemize}
  \item[-] Properties $1$ and $2$ in the Adjunction Conjecture hold,
            with the exception of precise inverse of adjunction, which
            is reduced to the divisorial case.
  \item[-] The moduli part $M_W$ is nef. Moreover, $M_W$ is semiample
           if $codim(j(W),X)=2$, according to \cite{Ka2}.
\end{itemize}

If $W$ is a curve, we just need Finite Base Change 
(Theorem~\ref{finite_basechange}) instead of the Base Change Conjecture.
Therefore the Adjunction Conjecture is proved if 
$dimX \le 3$, with the exception of precise inverse of adjunction 
(boundness was basically proven by Y. Kawamata \cite{Ka2}).

\section*{Acknowledgments}
This paper is part of author's dissertation submitted to
The Johns Hopkins University. I am indebted to my advisor, Professor
Vyacheslav V. Shokurov, who was so generous with his ideas, time 
and expertize. I am furtunate to have had him introduce me to this 
exciting area of research. 
%\newline
The pioneer work of Professor Yujiro 
Kawamata has motivated this project from the beginning. I wish to
thank him for updating me with his research.
\newline
I would also like to express my gratitude to those who have guided and
councelled me during my years of graduate education. Special
thanks go to Professors Bernard Shiffman and Steven Zucker.

\footnotetext[1]{This work was partially supported by NSF Grant 
DMS-9800807}

\newpage

%%%%%%%%%%%%%%%%%%%%%%%%%%%%%
%%%%%%%%%%%%%%%%%%%%%%%%%%%%%

\section{The basics of log pairs}   

%%%%%%%%%%%%%%%%%%%%%%%%%%%%%
%%%%%%%%%%%%%%%%%%%%%%%%%%%%%

%%%%%%%%%%%%%%%%%%%%%%%%%%%%%
%%%%%%%%%%%%%%%%%%%%%%%%%%%%%
\medskip
\subsection{Prerequisites}
  
%%%%%%%%%%%%%%%%%%%%%%%%%%%%%
%%%%%%%%%%%%%%%%%%%%%%%%%%%%%

A {\em variety} is a reduced irreducible scheme of finite type over
a fixed field $k$, of characteristic $0$.
\par
Let $X$ be a normal variety and $K$ one of the rings $\Z$, $\Q$ or $\R$.
A $K$-{\em divisor} $B=\sum_i b_i B_i$ on $X$
is a linear combination of prime Weil divisors with coefficients in $K$, 
i.e. an element of $N^1(X) \otimes K$. A $K$-divisor is said to be 
$K$-Cartier if it belongs to $Div(X) \otimes K \subset N^1(X) \otimes K$, 
where $Div(X)$ is the space of $\Z$-divisors which are Cartier. 
A $K$-divisor $B=\sum_i b_i B_i$ is {\em effective} if $b_i \ge 0$ for 
every $i$.
\par The fundamental invariant of $X$ is its
{\em canonical class} $K_X$. It is a $\Z$-Weil divisor, uniquely 
determined up to linear equivalence. In what follows, the choice of $K_X$
in its class is irrelevant.
\par Two $K$-divisors $D_1, D_2$ are $K$-linearly equivalent, denoted
by $D_1 \sim_K D_2$, if $D_1-D_2$ belongs to 
$P(X) \otimes K \subset N^1(X) \otimes K$, where $P(X)$ is the group of 
principal $\Z$-divisors associated to nonzero rational functions on $X$.
\par A morphism $f:X \to Y$ is called a {\em contraction} if the natural
morphism $\cO_Y \to f_*\cO_X$ is an isomorphism. Also, $f$ is called
an {\em extraction} if it is a proper birational morphism of normal 
varieties.
\par 
We say that a $K$-divisor $D$ is $K$-{\it linearly trivial over} $Y$, denoted
$D \sim_{K,f} 0$, if there exists a $K$-Cartier divisor 
$D'$ on $Y$ such that $D\sim_K f^*D'$. If $f$ is a contraction, then
the $K$-class of $D'$ is uniquely determined by $D$.
A $K$-Cartier divisor $D$ on $X$ is called $f$-{\em nef} if $D.C \ge 0$
for every proper curve $C$ such that $f(C)=\mbox{point}$.

%%%%%%%%%%%%%%%%%%%%%%%%%%%%%%%%%%%
%%%%%%%%%%%%%%%%%%%%%%%%%%%%%%%%%%%
\medskip
\subsection{Log varieties and pairs} 

%%%%%%%%%%%%%%%%%%%%%%%%%%%%%%%%%%%
%%%%%%%%%%%%%%%%%%%%%%%%%%%%%%%%%%%

The objects of the log-category are the singular counterpart of the 
smooth varieties with smooth boundary. They appear naturally 
in birational geometry.

\begin{defn} A {\em log pair\/} $(X,B)$
  is a normal variety $X$ equipped with an $\R$-Weil divisor $B$ such that 
  $K+B$ is $\R$-Cartier. We will equivalently say that $K+B$ is a 
  {\em log divisor}. A {\em log variety\/} is a log pair $(X,B)$ such that 
  $B$ is effective. 
 We call $B$ the {\em pseudo-boundary\/} of the log pair. 
\end{defn}

\begin{defn}
  \begin{enumerate}
    \item A log pair $(X,B)$ has {\em log nonsigular support\/}
  if $X$ is nonsingular and if $B=\sum b_i B_i$, then
  $\cup_{b_i \ne 0} B_i$ is a union of smooth divisors intersecting 
  transversely (in other words, it has {\it simple normal crossings}).
     \item A {\em log resolution \/} of a log pair $(X,B)$
 is an extraction $\mu:\tilde{X} \to X$ such that $\tilde{X}$ is
 nonsingular and $\Supp(\mu^{-1}(B)) \cup Exc(\mu)$ is a simple
 normal crossing divisor.
  \end{enumerate}
\end{defn}

One of the fundamental birational operations is the {\em pull back of log
divisors}. If $\mu:\tilde{X} \to X$ is an extraction and $K+B$ is a 
log divisor on $X$, there exists a unique 
log divisor $K_{\tilde{X}}+B^{\tilde{X}}$ on $\tilde{X}$ such that
 \begin{itemize}
   \item[i)] $B^{\tilde{X}}=\mu^{-1}B$ on $\tilde{X} \backslash Exc(\mu)$,
   \item[ii)] $\mu^*(K+B)= K_{\tilde{X}}+B^{\tilde{X}}$.
 \end{itemize}
The divisor $B^{\tilde{X}}$ is called the {\sl log codiscrepancy divisor} of 
$K+B$ on $\tilde{X}$, making $(\tilde{X},B^{\tilde{X}})$ a log pair 
which is identical to $(X,B)$ from the singularities point of view.
If $\mu:\tilde{X} \to (X,B)$ is a log resolution, then the log pair 
$(\tilde{X},B^{\tilde{X}})$ has log nonsingular support.
In the sequel, when we say that $\mu:(\tilde{X},\tilde{B}) \to (X,B)$
is a log resolution, it is understood that $\tilde{B}=B^{\tilde{X}}$.

\begin{example} Let $\mu:\tilde{X} \to X$ be the blow-up of a subvariety
$W$ of $X$ of codimension $c$, both nonsingular. Then 
$0^{\tilde{X}}=(1-c)E$ 
where $E$ is the exceptional divisor.
Therefore the log variety $X=(X,0)$ is ``similar'' to the log pair 
$(\tilde{X},(1-c)E)$. This illustrates the need for allowing 
the coefficients of the pseudo-boundary to take negative values.
\end{example}

%%%%%%%%%%%%%%%%%%%%%%%%%%%%%%%%%%%
%%%%%%%%%%%%%%%%%%%%%%%%%%%%%%%%%%%
\medskip
\subsection{Singularities and log discrepancies} 

%%%%%%%%%%%%%%%%%%%%%%%%%%%%%%%%%%%
%%%%%%%%%%%%%%%%%%%%%%%%%%%%%%%%%%%

The class of log canonical singularities can be described as the
largest class in which the LMMP seems to work, or as the smallest 
class containing Iitaka's log varieties which is closed under blow-ups.

\begin{defn} The log pair $(X,B)$ has {\em log canonical singularities} 
({\em lc} for short) if
there exists a log resolution $\mu:(\tilde{X},\tilde{B}) \to (X,B)$ such 
that all the coefficients of $\tilde{B}$ are at most $1$. 
\par
We say that $(X,B)$ has {\em Kawamata log terminal singularities} 
({\em klt} for short) if there exists a log resolution 
$\mu:(\tilde{X},\tilde{B}) \to (X,B)$ such that
the coefficients of $\tilde{B}$ are all less than $1$.
\end{defn}
It is easy to check that once $\tilde{B}$ has one of the above properties
on a log resolution, it has the same property on any log resolution.
In particular, a log pair $(X,B)$ with log nonsingular support is
log canonical (Kawamata log terminal) iff $B$ has coefficients at most $1$
(less than $1$). Note that both classes of singularities defined 
above have local nature.

\par
The singularities of log pairs are naturally described in terms
of {\em log discrepancies}. Discrepancies are invariants attributed
to Miles Reid who introduced them as means to control the canonical
class of variety under a birational base change. A normalized version
of discrepancies was also introduced in \cite{Shif}.

\begin{defn}
Let $(X,B)$ be a log pair. Let $E \subset Y \stackrel{\mu}{\to} X$ be a 
prime divisor on an {\em extraction} of $X$. 
The {\em log discrepancy} of $E$ with respect
to $K+B$ (or with respect to $(X,B)$), is defined as
$$
a_l(E;X,B)=1-e
$$
where $e$ is the coefficient of $E$ in the log codiscrepancy divisor $B^Y$.
By definition, $a_l(E;X,B)=1$ if $E$ is not in the support of $B^Y$.
The {\em center of $E$ on $X$} is $\mu(E)$, denoted by $c_X(E)$.
\newline
The log discrepancy $a_l(E;X,B)$ depends only on the discrete 
valuation defined by $E$ on $k(X)$, in particular independent
on the extraction $Y$ where $E$ appears as a divisor. 
\newline
In this paper we will write $a(E;X,B)$ or $a(E;B)$, dropping 
the index $l$ and even the variety $X$ from the notation. 
However, $a(E;B)$ should not be confused with the standard notation
in the literature for the {\em discrepancy} of $K+B$ in $E$, which 
is equal to $-1+a_l(E;X,B)$.
\end{defn}

\begin{rem} In the above notation, the 
log discrepancies for prime divisors on $Y$ are 
uniquely determined by the formula
$$
\mu^*(K_X+B)=K_Y+\sum_{E \subset Y} (1-a(E;X,B))E
$$
where the sum runs over all prime divisors of $Y$. 
\end{rem}

\begin{example}\label{ld_exmp}
Let $(X,B)$ be a log pair with log nonsingular support and
$E$ the exceptional divisor on the blow-up of the nonsingular
subvariety $Z \subset X$. Then 
 $$a(E;X,B)=\codim(Z,X)-j+\sum_{j \in J} a(E_j;X,B),$$
 where $J$ is the set of components of $B$ containing $Z$ and
 $j=|J| \le codim(Z,X)$. 
In particular, if $a(E_j;X,B) \ge 0$ for every $j \in J$, then
$$a(E;X,B) \ge \min_{j \in J} a(E_j;X,B) \ge 0$$
\end{example}

%%%%%%%%%%%%%%%%%%%%%%%%%%%%%%%%%%%%%%%%%
%%%%%%%%%%%%%%%%%%%%%%%%%%%%%%%%%%%%%%%%%
\medskip
\subsection{Log canonical centers}

%%%%%%%%%%%%%%%%%%%%%%%%%%%%%%%%%%%%%%%%%
%%%%%%%%%%%%%%%%%%%%%%%%%%%%%%%%%%%%%%%%%

 Let $(X,B)$ be a log pair and $x \in X$ a closed point. 
\begin{enumerate} 
   \item \cite{Ka1} A {\em log canonical center} ({\em lc center} for short) 
        of $(X,B)$ is a closed subvariety $W \subset X$ such that 
        $a(\eta_W;X,B) \le 0$. The minimal element of the set
        $$\{W'; x \in W', a(\eta_{W'};X,B) \le 0\}$$
        if it exists, is called the {\em minimal lc center} at 
        $x \in X$. If $\{x\}$ is not an lc center, then $(X,B)$
        is log canonical in a neighborhood of $x$, and moreover,
        the minimal lc center at $x$ exists if there is an 
        effective divisor $B^o \le B$ such that $K_X+B^o$
        is $\R$-Cartier and Kawamata log terminal in a neighborhood
        of $x$

   \item (V. Shokurov) An lc center $W$ is called {\em exceptional} if 
         $a(\eta_W;B)=0$ and on a log resolution 
         $\mu:(\tilde{X},\tilde{B}) \to (X,B)$ there exists a unique divisor
         $E$ such that $W=c_X(E)$ and $a(E;B) \le 0$ 
         (in particular, $a(E;B)=0$ if $dimW>0$).
         The definition does not depend on the choice of the log resolution.
         It is the generic equivalent of pure log terminality.
   \item We say that $(X,B)$ has a {\em normalized minimal lc center at $x$}
         if $a(x;B) \ge 0$ and on a log resolution 
         $\mu:(\tilde{X},\tilde{B}) \to (X,B)$ there exists a unique divisor
         $E$ such that $x \in c_X(E)$ and $a(E;B) \le 0$. In particular, 
         $a(E;B)=0$ and $W=c_X(E)$ is the minimal lc center at $x$. 
         The definition does 
         not depend on the choice of the log resolution. Moreover,
         there is an open neighborhood $U$ of $x$ such that
         $LCS(X,B)|_U=W|_U$ as schemes. In particular, $W$ is
         the only irreducible component of $LCS(X,B)$ passing through $x$.  
\end{enumerate}

If $W$ is an lc center for $(X,B)$, there might be several
prime divisors $E$ with $c_X(E)=W$ and $a(E;B)=0$. Such
divisors are called {\em lc places over} $W$ \cite{Ka1}. 
In fact, we have either infinitely many lc places over $W$, or 
exactly one.
The latter holds precisely when $W$ is an {\em exceptional lc center}.
\par
The unique place is realized as a divisor on an extraction of $X$, and
if $(E_1 \subset X_1 \to X)$ and $(E_2 \subset X_2 \to X)$
are two such realizations, then the induced birational
morphism $\tau:X_1 \cdots > X_2$ sends $E_1$ onto $E_2$ and
extends to an isomorphism in the generic point of $E_1$.
All codimension $1$ lc centers are exceptional (hopefully
this does not cause any confusion).

\begin{lem}\label{pert_lem}\cite{Ka1}(Perturbation Lemma) 
Let $K_X+B^o$ and $K_X+B$ be two log divisors
on $X$ such that $0 \le B^0 \le B$ and $K_X+B^o$ is Kawamata log terminal
in a neighborhood of $x$. If $(X,B)$ is log canonical at $x$, with
$W$ the minimal lc center at $x$, there exists an effective
$\Q$-Cartier divisor $D$ such that $K_X+B^o+(1-\epsilon)B+\epsilon D$
is log canonical with normalized minimal lc center $W$ at $x$, for
every $0<\epsilon <1$. Moreover, if $H$ is a $\Q$-free divisor
on $X$, we can assume $D \sim_\Q H$.
\end{lem}

%%%%%%%%%%%%%%%%%%%%%%%%%%%%%%%%%%%%%%%%%
%%%%%%%%%%%%%%%%%%%%%%%%%%%%%%%%%%%%%%%%%
\medskip
\subsection{The LCS locus}

%%%%%%%%%%%%%%%%%%%%%%%%%%%%%%%%%%%%%%%%%
%%%%%%%%%%%%%%%%%%%%%%%%%%%%%%%%%%%%%%%%%

Let $(X,B)$ be a log pair. The {\em locus of log canonical singularities}
of $K_X+B$ \cite[3.14]{Sh_flips} is the union of all lc centers:
$$
LCS(X,B)=\bigcup_{\textnormal{ W lc center}}  W
$$
The name is slightly confusing, in the sense that $W$ might be an lc center,
without $K_X+B$ being log canonical in $\eta_W$. A correct notation, 
proposed by J. Koll\'ar, is $Nklt(X,B)$: the abbreviation
for the locus where $(X,B)$ is not Kawamata log terminal. 
However, we will use Shokurov's notation since it better reflects its main 
use: to provide an induction step in higher dimensional algebraic 
geometry.

V. Shokurov also introduced a scheme structure on $LCS(X,B)$, defined as
follows. Let $\mu:(\tilde{X},\tilde{B}) \to (X,B)$ be a log resolution
and $P$ the truncation of $\tilde{B}$ to its components with
coefficients at least $1$. Then
$$
\cI(X,B)=\mu_*\cO_{\tilde{X}}(\ulcorner -P  \urcorner)
$$
is a coherent ideal sheaf on $X$, independent of the choice of
the log resolution. If $B$ is effective, then
$\cI(X,B) \simeq \mu_*\cO_{\tilde{X}}(\ulcorner -\tilde{B} \urcorner)$.
The ideal sheaf $\cI(X,B)$ defines a closed subscheme structure
on $LCS(X,B) \subset X$. 
\newline
The most general form of Kawamata-Viehweg vanishing theorem is:

\begin{thm}~\label{vanishing_LCS} (V. Shokurov) Let $L$ be a Cartier 
divisor on a log variety $(X,B)$ such that 
$$L \equiv K_X+B+H$$ where $H$ is a 
nef and big $\R$-divisor. Then $H^j(X,\cI(X,B)\otimes \cO(L))=0$
for every $j>0$. In particular, we have the following natural surjection
$$
H^0(X,L) \to H^0(LCS(X,B),L|_{LCS(X,B)}) \to 0
$$
\end{thm}

The following simple result plays a crucial role in the inverse of 
adjunction:

\begin{thm}(Connectedness Lemma \cite[5.7]{Sh_flips}, \cite[17.4]{Ko1})
Let $\pi:X \to S$ be a contraction and $K_X+B$ a log divisor on $X$
with the following properties:
\begin{enumerate}
   \item $-(K_X+B)$ is $\pi$-nef and $\pi$-big,
   \item the components of $B$ with negative coefficients are
         $\pi$-exceptional.
\end{enumerate}  
Then the induced map $LCS(X,B) \to S$ has connected fibers.
\end{thm}

One of its applications is the following result of J. Koll\'ar:

\begin{prop}\label{variation}\cite[Corollary 7.8]{Ko2}
Let $\{K_X+B_c\}_{c\in C}$ be an algebraic 
family of log divisors on $X$, parametrized by a smooth curve $C$.
For each closed point $x \in X$, the following subset of $C$ is
closed:
$$\{c \in C; x \in LCS(X,B_c) \} \subset C$$
\end{prop}

\newpage

%%%%%%%%%%%%%%%%%%%%%%%%%%%%%%
%%%%%%%%%%%%%%%%%%%%%%%%%%%%%%

\section{Minimal log discrepancies}

%%%%%%%%%%%%%%%%%%%%%%%%%%%%%%
%%%%%%%%%%%%%%%%%%%%%%%%%%%%%%

The minimal log discrepancy of a log pair $(X,B)$ in a closed subvariety 
$W \subset X$ is an invariant introduced by V. Shokurov. It can
be interpreted as the ``dimension'' of the singularity of $(X,B)$ 
in $W$, although it distinguishes log canonical singularities only.

\begin{defn} (V.V. Shokurov) For a log pair $(X,B)$ and a closed subset 
$W \subseteq X$  the following invariants are defined:
\begin{itemize}
   \item[-] $a(W;X,B)=\inf\{a(E;X,B); \emptyset \ne c_X(E) \subseteq W\}$
         is called the {\it minimal log discrepancy} of $(X,B)$ in $W$.
   \item[-] $a(\eta_W;X,B)=\inf\{a(E;X,B); c_X(E)=W\}$
         is called the {\it minimal log discrepancy} of $(X,B)$ in the
         generic point of $W$.
\end{itemize}
\end{defn}
We have $a(\eta_W;X,B) \ge a(W;X,B)$ and strict inequality holds 
in general. In fact, $$a(\eta_W;X,B)= a(W \cap U;U,B|_U)$$ for some
generic open subset $U \subseteq X$ intersecting $W$. We
abbreviate $a(X;X,B)$ and $a(\{x\};X,B)$ by $a(X;B)$ and $a(x;B)$,
respectively, where $x \in X$ is a closed point.
The following lemma shows that the minimal log discrepancy is
a well defined nonnegative real number if $(X,B)$ is log canonical 
in a neighborhood of $W$, and is equal to $-\infty$ otherwise.

\begin{lem}\label{ml} \cite[Proposition 17.1.1]{Ko1}
\begin{enumerate}
  \item If $(X,B)$ is not log canonical in a neighborhood of $W$,
then 
$$a(W;X,B)=-\infty$$
  \item Assume that $(X,B)$ is log canonical in a neighborhood of $W$.
Let $(\tilde{X},\tilde{B})$ be a log resolution of $(X,B)$ such that 
$F_W=\mu^{-1}(W)_{red}$ is a divisor and $\Supp(F_W)\cup \Supp(\tilde{B})$ 
has simple normal crossings. Then $a(W;X,B) \in \R_{\ge 0}$ and 
\begin{equation*}
\begin{split}
a(W;X,B)
& = \min\{a(F;X,B); F \textnormal{ irreducible component of } F_W \} \\
& = \sup\{c \ge 0; (\tilde{X},\tilde{B}+cF_W) 
                     \textnormal{ is log canonical near } W \} 
\end{split}
\end{equation*} 
Morever, the supremum is attained exactly on the components
of $F_W$ having log discrepancy minimal, that is equal to
$a(W;X,B)$.
\end{enumerate}
\end{lem} 

\begin{proof}
\begin{enumerate}
   \item If $(X,B)$ is not log canonical in a neighborhood of $W$,
          there exists a prime divisor $E$ on some extraction of $X$ 
          such that $a(E;B)<0$ and $W \cap c_X(E)  \ne \emptyset$. Let
          $x \in W \cap c_X(E)$ be a closed point. Suffices to show 
          that $a(x;B)=-\infty$, since $a(W;B) \le a(x;B)$. 
 Let $\mu:(\tilde{X},\tilde{B}) \to (X,B)$ 
 be a log resolution such that $E$ and $\mu^{-1}(x)$ are divisors on 
 $\tilde{X}$. Since $x \in \mu(E)$, there exists a component $E_0$ 
 of $\mu^{-1}(x)$ such that $E \cap E_0 \ne \emptyset$. Let
 $X_1$ be the blow up of $E \cap E_0$, with exceptional divisor $E_1$.
 Define inductively $X_k$ to be the blow up of $X_{k-1}$ in the intersection
 of $E$ and $E_{k-1}$, with exceptional divisor $E_k$. An easy 
 computation gives
 $$
 a(E_k;B)=k \cdot a(E;B)+a(E_0;B), \ \ c_X(E_k)=\{x\}
 $$
 for every $k \in \N$. Therefore 
 $\lim_{k \to \infty} a(E_k;B)=-\infty$, hence $a(x;B)=-\infty$.
  \item Shrinking $X$ to a neighborhood of $W$, we can assume that 
 $(X,B)$ is globally log canonical. 
 It is enough to check the invariance of the minimum under
 blow-ups on $\tilde{X}$, which follows from Example~\ref{ld_exmp}.
\end{enumerate}
\end{proof}

%%%%%%%%%%%%%%%%%%%%%%%%%%%%%%
%%%%%%%%%%%%%%%%%%%%%%%%%%%%%%
\medskip
\subsection{The lower semicontinuity of minimal log discrepancies}

%%%%%%%%%%%%%%%%%%%%%%%%%%%%%%
%%%%%%%%%%%%%%%%%%%%%%%%%%%%%%

Note what the first part of the proof of Lemma~\ref{ml} actually says:
$$
a(\eta_W;B)<0 \Longrightarrow a(x;B)<0 \mbox{ for every
closed point } x \in W
$$ 
This is typical for minimal log discrepancies: they are expected to 
behave in a lower semi-continuous fashion. 
To make this precise, fix a log pair $(X,B)$ and consider the function
$$
a:X \to \{-\infty\} \cup \R_{\ge 0},\ x \mapsto a(x;X,B)
$$

\begin{lem}
The nonempty set $\{x\in X;a(x) \ge 0\}$ is the biggest open subset of $X$
on which $(X,B)$ has log canonical singularities. Its closed complement
$\{x \in X; a(x)=-\infty\}$ is the union of all closed subvarieties
$W$ of $X$ such that $a(\eta_W;B)=-\infty$.
\end{lem}

\begin{proof} Let $x \in X$ such that $a(x) \ge 0$. According to the
observation above, $a(\eta_W;B) \ge 0$ for every closed subvariety of $X$
passing through $x$. Therefore there exists an open neighborhood $V$ of 
$x$ such that $(V,B|_V)$ is log canonical. In particular, $a(x') \ge 0$
for every $x' \in V$. Therefore $U=\{x\in X;a(x) \ge 0\}$ is open and
$(X,B)$ has log canonical singularities on $U$. The maximal property of 
$U$ is clear and the next lemma shows that
$U \ne \emptyset$.
\end{proof}

\begin{lem}\label{ns} We have $a(x)=dimX$ if 
$x \in Reg(X)\setminus \Supp(B)$ 
and $a(x)< dimX$ if $x \in Reg(X) \cap \Supp(B)$ and $B$ is
effective.
In particular, $a(x)$ is constant function equal to $dimX$ on an 
open dense subset of $X$.
\end{lem}

\begin{proof}
Let $x \in X$ a nonsingular point and let $E \subset \tilde{X} \to X$ 
be the exceptional divisor on the blow-up of $x$. Then 
$a(x;B) \le a(E;B)=dimX-\mult_x(B) \le dim X$, since $B$ is
effective.
\newline
If $x \in Reg(X)\setminus \Supp(B)$, then $a(x;B)=a(E;B)=dimX$
according to Lemma~\ref{ml}$.2$.
\end{proof}

The following conjecture gives the effective upper bound for the function
$a(x)$.

\begin{conj}~\label{max_mld}(V.V. Shokurov \cite{Sh_conj} ) 
Let $(X,B)$ be a log variety. Then 
$$
\sup_{x \in X} a(x) = dimX
$$ 
Moreover, the supremum is attained exactly on $Reg(X)\setminus \Supp(B)$. 
\end{conj}

The first part of Conjecture~\ref{max_mld} can be reduced to the following
conjecture.

\begin{conj}[Lower semi-continuity]\label{lower_semicontinuity} 
Let $(X,B)$ be a log variety.
Then the function $a(x)$ is lower semi-continuous, i.e. every closed 
point $x \in X$ has a neighborhood $x \in U \subseteq X$ such that 
$$
a(x;X,B)=\inf_{x' \in U} a(x';X,B)
$$
\end{conj}

Indeed, the function $a(x)$ may jump only downwards 
in special points and it is constant equal to $dimX$ on
an open dense subset of $X$. Therefore $\sup_{x \in X} a(x) = dimX$. 

\begin{example} We check Conjecture~\ref{lower_semicontinuity} for 
$dimX \le 2$. We may assume that
 $(X,B)$ is log canonical near $P \in X$, since otherwise there is nothing to 
 prove. In particular, $B$ is effective with coefficients at most $1$.
\begin{itemize}
   \item[a)] Assume $dimX=1$, and let $B=\sum_j (1-a_j)P_j$. 
             Then 
\[
a(x)= \left\{ \begin{array}{ll}
                   1    & \mbox{ if $x \notin \{P_j\}_j$} \\
                   a_j    & \mbox{ if $x =P_j$} 
                   \end{array} 
          \right. \]
            Note that the effectivity of $B$ is essential for lower 
semi-continuity.

   \item[b)] Assume $dimX=2$, and let $B=\sum_j (1-a_j)B_j$. We may assume 
             that $P \in \Supp(B_j)$ for every $j$,  
             $U=X\setminus P$ is nonsingular and $B|_U$ has
             nonsingular support. Then  
\[
a(x)= \left\{ \begin{array}{ll}
                   2    & \mbox{ if $x \in U \setminus \Supp(B)$} \\
                   1+a_j    & \mbox{ if $x \in U \cap B_j$}\\
              \end{array} 
          \right. \]
Now we have two cases. If $P\in X$ is a nonsingular point then
$a(P)\le \sum_j (a_j-1)+2 \le 1+a_j$ for every $j$. 
If $P\in X$ is singular, then it is well known that $a(P) \le 1$ and
equality holds iff $P \notin \Supp(B)$ and $(P\in X)$ is a DuVal 
singularity (cf. \cite[3.1.2]{Al2}).
   \end{itemize}
\end{example}

There are other interesting spectral properties of the minimal 
log discrepancies, conjectured by V. Shokurov, such as a.c.c. 
(see \cite{Sh_conj, Ko1} for details). Minimal log discrepancies were
successfully used by V. Shokurov for the existence and termination
of log flips in dimension $3$ \cite{Sh_flips}.

%%%%%%%%%%%%%%%%%%%%%%%%%%%%%%
%%%%%%%%%%%%%%%%%%%%%%%%%%%%%%
\medskip
\subsection{Precise inverse of adjunction}

%%%%%%%%%%%%%%%%%%%%%%%%%%%%%%
%%%%%%%%%%%%%%%%%%%%%%%%%%%%%%

\begin{conj}\label{precise} (cf. \cite[Conjecture 17.3]{Ko1})
Let $(X,B)$ be a log variety with a normalized
minimal lc center $W$ at $x$.
Let $\mu:(\tilde{X},\tilde{B}) \to (X,B)$ be a log 
resolution with $E$ the only lc center above $W$. 
Set $B_E=(\tilde{B}-E)|_E$. Then
$$
a((\mu|_E)^{-1}(x);E,B_E) = a(x;X,B).
$$

\[ \xymatrix{
   \tilde{X} \ar[d]_{\mu} & E \ar[d]^{\mu} \ar[l]_{\supset}\\
   X  & W \ar[l]^{\supset} 
} \]
\end{conj}

First of all, it is clear that $a((\mu|_E)^{-1}(x);E,B_E) \ge a(x;X,B)$,
so we just have to prove the opposite inequality:
$$
a(E \cap \mu^{-1}(x);E,B_E) \le  a(x;X,B)
$$
According to \cite[Chapter 17]{Ko1}, this inequality is implied by 
the Log Minimal Model Program in the case $0 \le a(x;X,B) \le 1$.

We describe in the next lemma a naive approach. Let $\tilde{B}=E+A$, 
hence $B_E=A|_E$. Assume $a(x;X,B) \ge 0$, $F_x=(\mu^{-1}(x))_{red}$ is a 
divisor and $\Supp(\tilde{B}) \cup \mu^{-1}(x)$ has simple normal crossings. 

\begin{lem}~\label{naive} With the above notations, assume moreover there 
exists an effective $\R$-divisor $\tilde{F}_x$ supported in $\mu^{-1}(x)$ 
with the following properties:
 \begin{enumerate}
   \item $-\tilde{F}_x$ is $\mu$-nef,
   \item the supremum 
         $a=\sup\{\alpha \ge 0; K_{\tilde{X}}+\tilde{B}+\alpha 
         \tilde{F}_x
         \mbox{ is log canonical above } x\}$
         is attained only on components $F$ with $a(F;B)=a(x;B)$.
 \end{enumerate}
Then $a((\mu|_E)^{-1}(x);E,B_E) \le a(x;X,B)$.
\end{lem}

\begin{proof}
Since $-(K_{\tilde{X}}+\tilde{B}+a \tilde{F}_x)=
          -\mu^*(K_X+B)-a \tilde{F}_x \equiv_{\mu} 
              - a \tilde{F}_x$
is $\mu$-nef, the induced map
$$
LCS(\tilde{X},\tilde{B}+a \tilde{F}_x) \to X
$$ 
has connected fibers from the Connectedness Lemma. The only candidates
for components of $LCS(\tilde{X},\tilde{B}+a  \tilde{F}_x)$ are 
$E$ and components of $F_x$ where $a$ is attained.
Therefore there exists a component $F$ of $F_x$ such that 
$F \cap E \ne \emptyset$ and $a(F;B)=a(x;B)$. Finally, 
$a(F \cap E;E,B_E)=a(F;\tilde{X},\tilde{B})=
a(x;X,B)$, hence
the desired inequality.
\end{proof}

Note that condition~\ref{naive}$.2$ is implied by the following
$$
(2') \ \ a(F_1;B)>a(F_2;B) \Longrightarrow 
\frac{a(F_1;B)}{\mult_{F_1}(\tilde{E}_x)} >
\frac{a(F_2;B)}{\mult_{F_2}(\tilde{E}_x)} 
$$
for any two components $F_1, F_2$ of $F_x$. In particular,
$F_x$ has this property. Unfortunately,
$-F_x$ is not $\mu$-nef in general.

\begin{example} If $x \in Reg(X)\setminus Supp(B)$ and $\mu$ is the blow 
up of $x \in X$, with exceptional divisor $E$, 
then $\tilde{F}_x=E$ satisfies the assumptions of the lemma. 
\end{example}

The following partial result on inverse of adjunction is due to 
V. V. Shokurov \cite[3.2]{Sh_flips} in the case $dimX=3$.
J\'anos Koll\'ar later found a formal proof based on the Connectedness
Lemma, which also proves the following result.

\begin{thm} (cf. \cite[Theorems 17.6,17.7]{Ko1}) Let $(X,B)$ be a log variety, 
$W$ an exceptional lc center with lc place $E$ and $x \in W$ a closed point.
\begin{enumerate}
   \item $W$ is the minimal lc center at $x$ for $(X,B)$ iff
         $(E,B_E)$ is Kawamata log terminal over a neighborhood 
         of $x$ in $X$.
   \item Assume that $W$ is the minimal lc center at $x$ for $(X,B)$
         and let $D$ be an effective $\R$-Cartier divisor on $X$ whose
         support does not contain $W$. Then
         $(X,B+D)$ is log canonical at $x$ iff $(E,(B+D)_E)$
         is log canonical above $x$.
\end{enumerate}
\end{thm}

\begin{proof}  
\begin{enumerate}
   \item Assume that $E$ is realized as a divisor on the log resolution 
         $\mu:(\tilde{X},\tilde{B}) \to (X,B)$, with $\tilde{B}=E+A$.
\[ \xymatrix{
   \tilde{X} \ar[d]_{\mu} & E \ar[d]^{\mu} \ar[l]_{\supset}\\
   X  & W \ar[l]^{\supset} 
} \]
         Note first the following equivalences:
\begin{itemize}
    \item[a)] $(E,B_E)$ is Kawamata log terminal over a neighborhood 
         of $x$ in $X$ iff $\mu^{-1}(x) \cap E \cap A^{\ge 1} =\emptyset$.
    \item[b)] $W$ is the minimal lc center at $x$ for $(X,B)$ iff
         $\mu^{-1}(x) \cap A^{\ge 1}=\emptyset$.
\end{itemize}
Therefore the implication ``$\Rightarrow$'' is clear. For the 
converse, assume that $\mu^{-1}(x) \cap E \cap A^{\ge 1}=\emptyset$.
Then $\mu^{-1}(x) \cap E \ne \emptyset$ and 
$\mu^{-1}(x) \cap A^{\ge 1}$ partition the set 
$\mu^{-1}(x) \cap (E \cup \Supp(A^{\ge 1}))$. 
\newline
But $E \cup \Supp(A^{\ge 1})=LCS(\tilde{X},\tilde{B})$ and
since $-(K_{\tilde{X}}+\tilde{B}) \equiv_\mu 0$ and
$\mu_*{\tilde{B}}=B$ is effective, the Connectedness Lemma
implies that $\mu^{-1}(x) \cap LCS(\tilde{X},\tilde{B})$ is a connected 
set. Therefore $\mu^{-1}(x) \cap A^{\ge 1}=\emptyset$. 

\item
The implication ``$\Rightarrow$'' is clear in $ii)$, so we only prove 
the converse. Assume that $(E,(B+D)_E)$ is log canonical above 
$x$. Shrinking $X$ near $x$, we can assume $(E,(B+D)_E)$
is globally log canonical and $(E,B_E)$ is Kawamata log terminal. 
Since $(B+D)_E=B_E+\mu^*(D|_W)$ and $D$ is effective, 
$(E,(B+tD)_E)$ is Kawamata log terminal for every $t<1$.
From $i)$, $W$ is the minimal lc center at $x$ for 
$(X,B+tD)$ for every $t<1$. In particular $a(x;B+tD) \ge 0$ 
for every $t<1$, hence for $t=1$ too.
\end{enumerate}
\end{proof}

\newpage

%%%%%%%%%%%%%%%%%%%%%%%%%%%%%%%%%%%%%%%%%
%%%%%%%%%%%%%%%%%%%%%%%%%%%%%%%%%%%%%%%%%

\section{Adjunction for fiber spaces}

%%%%%%%%%%%%%%%%%%%%%%%%%%%%%%%%%%%%%%%%%
%%%%%%%%%%%%%%%%%%%%%%%%%%%%%%%%%%%%%%%%%

Let $f:X \to Y$ be a proper contraction of normal varieties and
$K_X+B$ a log divisor which is log canonical over the generic
point of $Y$. We first introduce $B_Y$, the {\em discriminant} of 
the log divisor $K_X+B$ on $Y$. 
If $K_Y+B_Y$ is $\R$-Cartier, we could say that $K_Y+B_Y$ is 
the divisorial push-forward of the log divisor $K_X+B$.
\par
Restricting afterwards our attention to the case when 
$K_X+B \sim_{\R,f} 0$ and $K_Y+B_Y$ is a log divisor, we will
study
\begin{itemize}
  \item[a)] the relation between the singularities of $(X,B)$ and
$(Y,B_Y)$;
  \item[b)] the positivity properties of the $\R$-class
$M \in \Pic(Y) \otimes \R$ uniquely defined by the adjunction
formula
$$
K_X+B-f^*(K_Y+B_Y) \sim_\R f^*(M)
$$
\end{itemize}
We say that $K_Y+B_Y+M$ is the {\em push-forward} of
the log divisor $K_X+B$ on $Y$, being a combination of the
{\em divisorial part} $K_Y+B_Y$ and the {\em moduli part} $M$.
\par
Naturally, the push forward should be the inverse of the
pull back operation. The latter is naturally defined when
$f:X \to Y$ is a finite or birational morphism: if 
$K_Y+B$ is a log divisor on $Y$, 
there exists an induced log divisor $K_X+B^X$ on $X$ uniquely
defined by the adjunction formula
$$
K_X+B^X=f^*(K_Y+B)
$$ 
Then $(B^X)_Y=B$ and $M=0$, that is $B$ is the discriminant of $K_X+B^X$
on $Y$ and the moduli part is trivial.

%%%%%%%%%%%%%%%%%%%%%%%%%%%%%%%%%%%%%%%%%
%%%%%%%%%%%%%%%%%%%%%%%%%%%%%%%%%%%%%%%%%
\medskip
\subsection{The discriminant of a log divisor}

%%%%%%%%%%%%%%%%%%%%%%%%%%%%%%%%%%%%%%%%%
%%%%%%%%%%%%%%%%%%%%%%%%%%%%%%%%%%%%%%%%%

The following is the invariant form of the definition proposed by 
Y. Kawamata in \cite{Ka2,Ka3}.

\begin{defn}\label{defnpf} Let $f:X \to Y$ be a surjective 
morphism of normal
varieties and $K_X+B$ a log divisor which is log canonical over
the generic point of $Y$. For a prime divisor $Q \subset Y$ define
$$
a_Q=\sup\{c \in \R; K_X+B+c f^*Q \mbox{ is log canonical over } 
      \eta_Q\}
$$
Then $B_Y=\sum_Q (1-a_Q) Q$ is a well defined $\R$-Weil divisor 
on $Y$, called the {\em discriminant} of the log divisor $K_X+B$ on 
$Y$. 
\end{defn}

\begin{rem}\label{basic}
\begin{enumerate}
  \item By abuse of language, $f^*Q$ is defined as the divisor 
 associated to the pullback $f^*t$ of a local parameter $t$ of $Q$ on $Y$. 
 Since the supremum is defined over the generic point of $Q$, the choice of 
 $t$ is irrelevant.
   \item If $f':(X',B'=B^{X'}) \to (X,B) \to Y$ is the map
 induced by a crepant extraction or a finite cover of $(X,B)$, 
 then $B'_Y=B_Y$. In other words, for computing $B_Y$ we are free to 
 replace $(X,B)$ by any crepant extraction or finite cover $(X',B^{X'})$.
   \item $B_Y$ is well defined since $a_Q=1$ for all but a finite
number of prime divisors. Indeed, assuming that $X$ is nonsingular,
there exists a non-empty open subset $U \subseteq Y$ such that
$K_X+B$ is log canonical over $U$ and $f$ has nonsingular 
(possibly disconnected) fibers over $U$. This implies $a_Q=1$ for every
prime divisor $Q$ with $Q \cap U \ne \emptyset$.
   \item For any prime divisor $Q$, $a_Q$ is a real number because 
$K_X+B$ is log canonical over $\eta_Y$.
To compute $a_Q$, we can assume after blowing up $X$ that
$(X,B_{red}+(f^*Q)_{red})$ is log nonsingular over an open subset 
$U \subseteq Y$ with $U \cap Q \ne \emptyset$. Let $f^*Q=\sum_j w_j P_j$ 
over $\eta_Q$ (note that $f(P_j)=Q$ for every $P_j$). Then
$$
a_Q=\min_j \frac{a(P_j;B)}{w_j}
$$
In other words, if $b_j=\mult_{P_j}(B)$ for every $j$, then
$$
b_Q=1-a_Q=\max_j \frac{b_j+w_j-1}{w_j}
$$
This is exactly the formula proposed in \cite{Ka2, Ka3}. 
In particular, $B_Y$ has rational coefficients if $B$ does. 
  \item In the above notation,
$
 \frac{1}{N} \min_j a(P_j;B)  \le a_Q \le \min_j a(P_j;B),
$
where $N=\max_j w_j \in \N$. The presence of $N$ makes the precise
inverse of adjunction for higher codimension lc centers an inequality
instead of equality. 
  \item(Additivity) If $D$ is an $\R$-Cartier divisor on $Y$, 
 then $K_X+B+f^*D$ is again log canonical over $\eta_Y$ and 
 $(B+f^*D)_Y=B_Y+D$. 
\end{enumerate}
\end{rem}

\begin{example}\label{exmp_pf}
\begin{enumerate}
 
  \item Assume that $f$ is birational and $K_X+B \sim_{\R,f} 0$. 
        Then $B_Y=f_*(B)$, $K_Y+B_Y$ is $\R$-Cartier and 
        $f:(X,B) \to (Y,B_Y)$ becomes a crepant extraction, that is 
        $K_X+B=f^*(K_Y+B_Y)$.
 
  \item Assume that $f$ is a finite map and $K_X+B^X$ is the pull back
        of the log divisor $K_Y+B$. Then 
             $$(B^X)_Y=B$$ Moreover, if 
        $P$ is a prime divisor on $X$, $Q=f(P)$ and $w=\mult_P(f^*Q)$,
        then $a_Q=a(Q;Y,B)=a(P;X,B^X)/w$.
  
  \item Assume $f:X \to Y$ is a fiber space of smooth varieties
        with simple normal crossing ramification. Then $f$ is
        semistable in codimension $1$ iff $K_X=K_X+0$ has
        trivial discriminant on $Y$.

 \item Assume that $Y=C$ is a smooth curve. Then $f:(X,B) \to C$ is 
       {\it log canonical} in the sense of \cite[Definition 7.1]{KM}
       iff $(C,B_C)$ has canonical singularities, that is $B_C \le 0$.
\end{enumerate}
\end{example}

The following result of Y. Kawamata gives a cohomological 
sufficient condition for the effectivity of the discriminant.

\begin{lem}\label{effectivity} \cite{Ka3} Let $f:X \to Y$ be a surjective
            map and $K_X+B$ a log divisor with log nonsingular support
            which is Kawamata log terminal over $\eta_Y$. 
            Let $Q$ be prime divisor on $Y$ such that 
            the coefficient of $Q$ in the discriminant $B_Y$ is negative. 
            Then the following hold:
  \begin{itemize}
      \item[-] $\ulcorner -B \urcorner$ is effective over
            the generic point of $Q$.
      \item[-] The induced map
            $\cO_{Y,Q} \to (f_*\cO_X(\ulcorner -B \urcorner))_Q$
            is not surjective.
  \end{itemize}
\end{lem}

\begin{proof}
Since $\mu_*\cO_{\tilde{X}}(\ulcorner -\tilde{B} \urcorner)=
\cO_X(\ulcorner - B \urcorner)$ for any log resolution 
$\mu:(\tilde{X},\tilde{B}) \to (X,B)$, we can assume that
we are in the situation and notation of Remark~\ref{basic}$.4$.
Then $b_Q<0$ is equivalent to 
$$\ulcorner - b_j \urcorner \ge w_j \mbox{ for every } j$$
Since $\ulcorner - B^h \urcorner \ge 0$, there exists 
an open subset $U \subset X$ such that $U \cap Q \ne \emptyset$
and $\ulcorner -B \urcorner$ is effective on $\mu^{-1}(U)$.
Moreover, the induced inclusion 
$\cO_{Y,Q} \to (f_*\cO_X(\ulcorner -B \urcorner))_Q$
factors as:
$$
{\cO}_{Y,Q} \to {{\cO}_Y(Q)}_Q \to 
        f_*{\cO}_X(\ulcorner - B \urcorner)_Q
$$
In particular, the natural inclusion 
$\cO_{Y,Q} \to (f_*\cO_X(\ulcorner -B \urcorner))_Q$
cannot be surjective.
\end{proof}

%%%%%%%%%%%%%%%%%%%%%%%%%%%%%%%%%%%%%%%%%
%%%%%%%%%%%%%%%%%%%%%%%%%%%%%%%%%%%%%%%%%
\medskip
\subsection{Base change for the divisorial push forward}

%%%%%%%%%%%%%%%%%%%%%%%%%%%%%%%%%%%%%%%%%
%%%%%%%%%%%%%%%%%%%%%%%%%%%%%%%%%%%%%%%%%

The following result shows that the divisorial push forward of
a log divisor commutes with finite base changes.

\begin{thm}[Finite Base Change] ~\label{finite_basechange} 
Let $f:X \to Y$ be a proper
surjective morphism of normal varieties and $K_X+B$ a $K$-Cartier log
divisor which is log canonical over $\eta_Y$. Let
$f:X' \to Y'$ be a morphism induced by the finite base change
$\sigma:Y' \to Y$, and set $B'=B^{X'}$. 
\[ \xymatrix{
   X \ar[d]_f & X' \ar[l]_{\nu} \ar[d]^{f'}\\
   Y          & Y' \ar[l]^{\sigma}
} \]
Then $K_{X'}+B'$ is log canonical over $\eta_Y$ and
$\sigma^*(K_Y+B_Y)=K_{Y'}+(B')_{Y'}$ as $K$-Weil divisors. 
Moreover, $K_Y+B_Y$ is $K$-Cartier iff 
$K_{Y'}+(B')_{Y'}$ is $K$-Cartier. 
\end{thm}

\begin{proof} To check the equality $\sigma^*(K_Y+B_Y)=K_{Y'}+(B')_{Y'}$
we choose an arbitrary prime divisor $Q' \subset Y'$, with $\sigma(Q')=Q$ 
and $\mult_{Q'}(\sigma^*Q)=w$. From Example $~\ref{exmp_pf}.2$, 
we have to show that $a_{Q'}=w \cdot a_Q$.
\newline
If $c \le a_Q$, then $K_X+B_X+cf^*Q$ is log canonical over $\eta_Q$, hence 
$K_{X'}+B'+c(f \circ \nu)^*Q=K_{X'}+B'+c {f'}^*\sigma^*Q$ is log canonical 
over $\eta_Q$. But $K_{X'}+B'+c {f'}^*\sigma^*Q \ge K_{X'}+B'+c w {f'}^*Q'$,
so $K_{X'}+B'+c w {f'}^*Q'$ is log canonical. Hence $cw \le a_{Q'}$.
In particular, $a_{Q'} \ge w \cdot a_Q$.
\newline
Conversely, let $c\ge a_Q$. After possible blow-ups on $X$,
there exists a prime divisor $P$ on $X$ with $a(P;B+cf^*Q)\le 0$ and
$f(P)=Q$. Since $X'$ is a resolution
of $X \times_Y Y'$, there exists a prime divisor $P'$
on $X'$ with $\nu(P')=P, f'(P')=Q'$. By Example $~\ref{exmp_pf}.2$, 
$a(P';B'+c \cdot w {f'}^*Q')=a(P';B'+c(f\circ \sigma)^*Q)\le 0$,
so $c \cdot w \ge a_{Q'}$. Therefore $w \cdot a_Q \ge a_{Q'}$.
\par
The rest follows from the next lemma.
\end{proof}
 
\begin{lem}
Let $\sigma:Y' \to Y$ be a finite map of
normal varieties, $D$ a $\Q$-Weil divisor on $Y$,
$D'=\sigma^*(D)$ the pull back of $D$, which is a $\Q$-Weil 
divisor on $Y'$, and $r \in \N$. Then
\begin{enumerate}
  \item if $rD$ is Cartier, then $rD'$ is Cartier;
  \item if $rD'$ is Cartier then $(deg(\sigma)\cdot r) D$ is Cartier.
\end{enumerate}
\end{lem}

Under certain conditions, we expect that the divisorial 
push forward commutes with birational base changes too. 
According to \cite[5.12,9.12]{Mo} and \cite{Ka2, Ka3}, we anticipate  
the following conjecture to be true. A partial result in this
direction is Proposition~\ref{sigma_prop}.

\begin{conj}[The Base Change Conjecture]~\label{base_change_conj} 
Let $f:X \to Y$ 
be a contraction of normal varieties and let $K_X+B$ be log divisor 
with the following properties:
  \begin{itemize}
    \item[-] $K_X+B \sim_{\Q,f} 0$;
    \item[-] $K_X+B$ is Kawamata log terminal over $\eta_Y$;
    \item[-] $(X,B)$ has log nonsingular support and 
             $\cO_{Y,\eta_Y}=
                 (f_*\cO_X(\ulcorner -B \urcorner))_{\eta_Y}$.
  \end{itemize}
Then $K_Y+B_Y$ is $\Q$-Cartier and if $f':(X',B^{X'}) \to Y'$ is 
             a contraction induced by a birational base change 
             $\sigma:Y' \to Y$, then $(B_Y)^{Y'}=(B^{X'})_{Y'}$. In
             other words,
$$
 \sigma^*(K_Y+B_Y)=K_{Y'}+(B^{X'})_{Y'}.
$$
\end{conj}

The divisor $K_Y+B_Y$ is always $\Q$-Cartier if $Y$ is 
$\Q$-factorial, in particular nonsingular. As for the base change,
even if it does not hold for $f:(X,B) \to Y$, it should hold
for data $f':(X',B') \to Y'$ induced on ``sufficiently large
extractions'' of $Y$. 
\par
 The Base Change Conjecture is intuitively equivalent to the 
Inverse of Adjunction Conjecture. As the next result 
shows, the log divisor $K_X+B$ and its divisorial push forward log 
divisor should be {\em in the same class of singularities}.

\begin{prop}\label{inv} Assume the Base Change Conjecture holds true
for $f:(X,B) \to Y$ and let $Z$ be a closed proper subset of $Y$. 
There exists a positive natural number $N \in \N$ such that
$$
\frac{1}{N}a(f^{-1}(Z);X,B) \le a(Z;Y,B_Y) \le 
a(f^{-1}(Z);X,B).
$$
\end{prop}

\begin{proof}
 There exists a fiber space $f':(X',B') \to Y'$ induced
by a birational base change $\sigma:Y' \to Y$ with the following properties:
 \begin{itemize}
   \item[-] $Y'$ is nonsingular, $\sigma^{-1}(Z)$ is a divisor and
            $\Supp(\sigma^{-1}(Z)) \cup \Supp((B_Y)^{Y'})$ is included
            in a snc divisor $Q=\sum_l Q_l$;
   \item[-] $X'$ is nonsingular and $\Supp(B') \cup \Supp({f'}^*Q)$ is
            included in a snc divisor $P=\sum_j P_j$;
   \item[-] there exists an index $j_0$ such that 
            $P_{j_0} \subseteq (\sigma \circ f')^{-1}(Z)$,
            $f'(P_{j_0})=Q_{l_0}$ for some index $l_0$ and
            $a(P_{j_0};X,B)=a(f^{-1}(Z);X,B)$.
 \end{itemize}

\[ \xymatrix{
   X \ar[d]_f & X' \ar[l]_{\nu} \ar[d]^{f'}\\
   Y          & Y' \ar[l]^{\sigma}
} \]
 
 Indeed, the first two property are obtained by letting
 $\sigma:Y' \to Y$, and then $X' \to X \times_Y Y'$ be 
 ``large enough'' resolutions. As for the third,
 let $P_{j_0}$ included in ${f'}^{-1}\sigma^{-1}(Z)$ such that 
 $a(P_{j_0};X,B)=a(f^{-1}(Z);X,B)$. If $f'(P_{j_0})$ is a divisor, there
 is some $l_0$ with $f'(P_{j_0})=Q_{l_0}$, so we are done. 
 Otherwise, by further blow-ups on $Y'$ and $X'$ we can assume the 
 proper transform of $P_{l'}$ maps to a divisor. 
 Note that we do not change any horizontal component, since we only 
 perform operations over proper subsets of $Y$.
 \par
 Then $a(Q_{l_0};B_Y)=a(Q_{l_0};B_{Y'}) \le a(P_{j_0};B')=a(f^{-1}(Z);B)$, 
 hence $$a(Z;B_Y) \le a(f^{-1}(Z);B)$$ 
 On the other hand, if $N=\max\{w_{lj}\} \in \N$ and $Q_l$ is any divisor 
 of $Q$ contained in $\sigma^{-1}(Z)$, then 
 $a(Q_l;B_Y)=a(Q_l;B_{Y'}) \ge \frac{1}{N} 
  \min_{f'(P_j)=Q_l} a(P_j;B') \ge \frac{1}{N} a(f^{-1}(Z);B)$.
 Taking infimum after all these $Q_l$'s we obtain the other inequality.
\end{proof}

%%%%%%%%%%%%%%%%%%%%%%%%%%%%%%%%%%%%%%%%%
%%%%%%%%%%%%%%%%%%%%%%%%%%%%%%%%%%%%%%%%%%
\medskip
\subsection{Positivity of the moduli part}

%%%%%%%%%%%%%%%%%%%%%%%%%%%%%%%%%%%%%%%%%%
%%%%%%%%%%%%%%%%%%%%%%%%%%%%%%%%%%%%%%%%%%

Let $f:(X,B) \to Y$ be a data satisying the assumptions of the
Base Change Conjecture. Assuming that $K_Y+B_Y$ is $\Q$-Cartier,
there exists a unique class $M_Y \in \Pic(Y) \otimes \Q$ satisying
the following adjunction formula:
$$
K_X+B \sim_\Q f^*(K_Y+B_Y+M_Y).
$$
We can rewrite the above formula as 
$K_{(X,B)/(Y,B_Y)}:=K_X+B-f^*(K_Y+B_Y) \sim_\Q f^*M$.
Thinking of $f:(X,B) \to (Y,B_Y)$ as being the log analogue of semistable in 
codimension $1$ morphisms, the line bundle $\cO_Y(\nu M)$ corresponds to 
$f_*\cO_X(\nu K_{(X,B)/(Y,B_Y)})$ for divisible and large enough 
integers $\nu \in \N$. Therefore we expect the following conjecture
on the positivity of log-Hodge bundles to be true.

\begin{conj}[Positivity]\label{pos_conj}
(cf. \cite{Ka2, Ka3},\cite[5.12,9.12]{Mo}) 
Let $f:X \to Y$ be a contraction of 
normal varieties and let $K_X+B$ be log divisor with the following 
properties:
  \begin{itemize}
    \item[-] $K_X+B \sim_{\Q,f} 0$;
    \item[-] $K_X+B$ is Kawamata log terminal over $\eta_Y$;
    \item[-] $(X,B)$ has log nonsingular support and 
             $\cO_{Y,\eta_Y}=
                 (f_*\cO_X(\ulcorner -B \urcorner))_{\eta_Y}$.
  \end{itemize}
Then $K_Y+B_Y$ is $\Q$-Cartier and the moduli part $M_Y$ is $\Q$-free,
that is the line bundle $\cO_Y(\nu M_Y)$ is generated by global sections
for some integer $\nu \in \N$.
\end{conj}

\begin{rem} The connectivity assumption on the fibers of $f$ 
is essential for the positivity of $M$. For instance, assume 
$f:X \to Y$ is a finite morphism between two smooth projective 
curves and let $K_X+B$ be a log divisor on $X$.
Then $K_Y+B_Y$ is a log divisor and 
$$K_X+B_X-f^*(K_Y+B_Y)=-\Delta$$ where $\Delta$ is an effective
divisor whose support does not contain any set theoretic
fiber of $f$. Therefore $-M$ is nef.
An important particular case is when $f$ is a Galois cover and
$K_X+B$ is Galois invariant. Then $M=0$.
\end{rem}

Under an extra assumption, Y. Kawamata proved that the moduli part 
$M_Y$ is a nef divisor \cite{Ka3}, and moreover, $M_Y$ is $\Q$-free if 
$f$ has relative dimension $1$ \cite{Ka2}.
We end this section with his positivity result.

\begin{thm}\label{technical} (cf. \cite[Theorem 2]{Ka3})
Consider the following setting:
\begin{enumerate}
  \item $f:X \to Y$ is a contraction of nonsingular projective varieties;
  
  \item $K_X+B \sim_{\Q, f} 0$;

  \item there exist simple normal crossing divisors $P=\sum P_j$
        and $Q=\sum Q_l$ on $X$ and $Y$, respectively, such that
        $\Supp(B) \subset P$, $f^{-1}(Q) \subset P$ and $f$ is smooth 
        over $Y \backslash Q$;

  \item $B = B^h + B^v$ such that any irreducible component of $B^h$ is 
        mapped surjectively onto $Y$ by $f$, $f: \text{Supp}(B^h) \to Y$ 
        is relatively normal crossing over $Y \setminus Q$, and
        $f(\text{Supp}(B^v)) \subset Q$.  An irreducible 
        component of $B^h$ (resp. $B^v$) is called {\it horizontal} 
        (resp. {\it vertical});

  \item $K_X+B$ is Kawamata log terminal over $\eta_Y$ and
        $\cO_{Y,\eta_Y}= 
        (f_*{\mathcal O}_X(\ulcorner - B \urcorner))_{\eta_Y}$. 
\end{enumerate}
Then the moduli part $M_Y$ is a nef divisor on $Y$.
\end{thm}

\begin{rem} In \cite[Theorem 2]{Ka3} it is further assumed 
that $\ulcorner -B \urcorner$ is effective (that is $(X,B)$
has Kawamata log terminal singularities), although this assumption
is not used in the proof. Indeed, set $D=B-f^*B_Y$ and
$$
Z=\bigcup \{f(P_j);P_j \subset \textnormal{Supp}(B), 
             \text{codim}(f(P_j),Y) \ge 2\}.
$$ 
The following hold:
\begin{itemize}
    \item[-] $K_X+D \sim_\Q f^*(K_Y+M_Y)$;
    \item[-] $D_Y=0$;
    \item[-] $\ulcorner -D \urcorner$ is effective over $Y \backslash Z$;
    \item[-] $Z \subset Q$ and $\textnormal{codim}(Z,Y) \ge 2$.
\end{itemize}
Including $Z$ in the closed subset of codimension at least $2$ that
is disregarded throughout the proof of \cite[Theorem 2]{Ka3}, we 
obtain the nefness of $M_Y$.
\end{rem}

We say that $f:(X,B) \to Y$ has the property $(\star)$ if the assumptions
of the above theorem hold true.

\begin{prop}\label{sigma_prop} Let $f:(X,B) \to Y$ and $f':(X',B') \to Y'$ 
be fiber spaces satisfying the property $(\star)$ such that the following 
hold:
\begin{enumerate}
  \item $f'$ is induced by $f$ by the birational base change $\sigma$,
        $\sigma^{-1}(Q) \subset Q'$ and $\nu^{-1}(P) \subset P'$;
  \item $B'=B^{X'}$.
\end{enumerate}
\[ \xymatrix{
   (X,B) \ar[d]_f & (X',B') \ar[l]_{\nu} \ar[d]^{f'}\\
   Y          & Y' \ar[l]^{\sigma}
} \]
Then $\Sigma=(B^{X'})_{Y'}- (B_Y)^{Y'}$ is an effective $\sigma$-exceptional 
divisor.
\end{prop}

\begin{proof}
 We have $\nu^*(K_X+B)=K_{X'}+B' \sim_\Q {f'}^*(K_{Y'}+(B')_{Y'}+M')$. 
On the other hand, $\nu^*(K_X+B) \sim_\Q \nu^* f^*(K_Y+B_Y+M)
={f'}^*{\sigma}^*(K_Y+B_Y+M)$.
Since $f'$ is a contraction, we infer that
$$
\sigma^*(K_Y+B_Y+M) \sim_\Q K_{Y'}+(B')_{Y'}+M'
$$
Therefore $\Sigma \sim_\Q -M'+\sigma^*M$. Clearly $\Sigma$ is 
$\sigma$-exceptional. Moreover, since $-\Sigma$ is $\sigma$-nef,
the negativity of the birational contraction $\sigma$ implies that 
$\Sigma$ is effective.
\end{proof}

\newpage

%%%%%%%%%%%%%%%%%%%%%%%%%%%%%%%%%%%%%%%%%%
%%%%%%%%%%%%%%%%%%%%%%%%%%%%%%%%%%%%%%%%%%

\section{Adjunction on log canonical centers}

%%%%%%%%%%%%%%%%%%%%%%%%%%%%%%%%%%%%%%%%%%
%%%%%%%%%%%%%%%%%%%%%%%%%%%%%%%%%%%%%%%%%%

Let $K+B$ be a log divisor on a normal variety $X$.
For a closed subvariety $W \subset X$ such that
$a(\eta_W;X,B)=0$, the Adjunction Conjecture predicts
that the different $B_{W^\nu}$ induced on the
normalization of $W$ has the following properties:
\begin{itemize}
   \item[-] $K_{W^\nu}+B_{W^\nu}$ is a log divisor with singularities 
             similar to those of $K+B$ near $W$;
   \item[-] The moduli part $M_{W^\nu}$, uniquely defined by the 
            adjunction formula
$$
(K_X+B)|_{W^\nu} \sim_\R K_{W^\nu}+B_{W^\nu}+M_{W^\nu}
$$
is an $\R$-free divisor.
\end{itemize}
This chapter contains partial results towards this conjecture.

%%%%%%%%%%%%%%%%%%%%%%%%%%%%%%%%%%%%%%%%%%
%%%%%%%%%%%%%%%%%%%%%%%%%%%%%%%%%%%%%%%%%%
\medskip
\subsection{The different}

%%%%%%%%%%%%%%%%%%%%%%%%%%%%%%%%%%%%%%%%%%
%%%%%%%%%%%%%%%%%%%%%%%%%%%%%%%%%%%%%%%%%%

\begin{defn} We say that $j:Y \to (X,B)$ is an {\em adjunction setting} 
if the following hold:
\begin{enumerate}
   \item $j:Y \to X$ is a proper morphism of normal varieties, generically
         one-to-one onto its image $j(Y)=W$;
   \item $K_X+B$ is a $K$-Cartier log divisor such that $a(\eta_W;X,B)=0$.
\end{enumerate}

\end{defn}

\begin{rem}
If $W^\nu$ is the normalization of $W$, then
$\nu:W^\nu \to (X,B)$ is an adjunction setting and there exists a unique
birational contraction $\sigma:Y \to W^\nu$ making
the following diagram commutative:
\[ \xymatrix{
       Y \ar[rr]^{\sigma} \ar[dr]_j &  & {W^\nu} \ar[dl]^\nu \\
         &                    (X,B)  &
} \]
Conversely, any such birational contraction $\sigma: Y \to W^\nu$
induces the adjunction setting  $\nu \circ \sigma:W^\nu \to (X,B)$.
\end{rem}

The main property of adjunction settings is that the log divisor 
$K_X+B$ has a natural {\em different} $B_Y\in N^1(Y) \otimes K$, measuring
the singularities of $(X,B)$ over the codimension $1$ points of $Y$.
To define the different, we assume that $j(Y)$ is an exceptional lc center 
of $K_X+B$. See Remark~\ref{semin} for the general case.

\begin{defn}\label{res_def}
Assume $j:Y \to (X,B)$ is an adjunction setting and $j(Y)$ is an exceptional
lc center. Let $\mu:(\tilde{X},\tilde{B}) \to (X,B)$
be a log resolution such that $E$, the lc place over $j(Y)$, is
realized as a divisor on $\tilde{X}$ and the induced map $E \to j(Y)$
factors through $Y$:
\[ \xymatrix{
   (E,B_E) \ar[d]_f \ar[r]^{\subset}& (\tilde{X},\tilde{B}) \ar[d]^{\mu}\\
   Y \ar[r]^j         & (X,B) 
} \]
Set $B_E=A|_E$, where $A=\tilde{B}-E$ is a divisor not containing $E$ in its 
support. Then the {\em different} of the log divisor $K_X+B$ on $Y$,
denoted $B_Y$, is defined as the discriminant on $Y$ of the log divisor 
$K_E+B_E$.
\end{defn}

\begin{rem}
   The definition is independent on the choice of the log resolution
$\mu$. Indeed, using Hironaka's hut, assume 
$\mu'=\mu \circ \tau:(\tilde{X'},\tilde{B'}) \to (X,B)$ 
is another log resolution induced by the extraction 
$\tau:\tilde{X'} \to \tilde{X}$. Let $E'$ and $E$ be the lc places above 
$j(Y)$ on $\tilde{X'}$ and $\tilde{X}$ respectively. By uniqueness,
$E'$ is the proper transform of $E$ via $\tau$ and
$(\tau|_{E'}):E' \to E$ is an extraction. In particular,
$\tau^*(K_{\tilde{X}}+E+A)=K_{\tilde{X'}}+E'+A'$, so the classical adjunction
formula gives
$$
(\tau|_{E'})^*(K_E+B_E)=K_{E'}+B_{E'}
$$
Therefore $(E,B_E)$ and $(E',B_{E'})$ have the same discriminant on $Y$ 
according to Remark~\ref{basic}$.2$.
\end{rem}

\begin{lem}\label{stan} In the notations of Definition~\ref{res_def}, assume
moreover that $f(Y)$ is an irreducible component of $LCS(X,B)$. Then
the following hold:
\begin{enumerate}
    \item $K_E+B_E \sim_K f^*j^*(K_X+B)$. 
    \item $f$ is a contraction.
    \item $(E,B_E)$ is Kawamata log terminal over $\eta_Y$ and 
          $\cO_{Y,\eta_Y} = (f_*\cO_E(\ulcorner -B_E \urcorner))_{\eta_Y}$.
    \item $a(Z;X,B) \le a((\mu|_E)^{-1}(Z);E,B_E)$ for every closed subset
           $Z \subseteq j(Y)$. 
\end{enumerate}
\end{lem}

\begin{proof}
\begin{enumerate}
    \item We have $K_E+B_E \sim_K (K_{\tilde{X}}+\tilde{B})|_E
               \sim_K \mu^*(K_X+B)|_E =f^*j^*(K_X+B)$.
    \item Since $W=j(Y)$ is an irreducible component of $LCS(X,B)$, 
          there exists an open subset $U \subseteq X$ such that
          $W \cap U$ is the only lc center for $(U,B|_U)$. Since
          $W$ is an exceptional center, 
          $LCS({\tilde{X}},\tilde{B})|_{\mu^{-1}(U)}=E|_{\mu^{-1}(U)}$,
          so the the Connectivity Lemma implies that $E \to W$ has
          connected fibers over $U\cap W \ne \emptyset$. 
          The induced morphism $E \to W^\nu$ has thus connected fibers,
          so the same holds for $f:E \to Y$.
    \item If $U=X \setminus \mu(\Supp(A^{\ge 1}))$ then 
          $W \cap U \ne \emptyset$ and $(E,B_E)$ is Kawamata log terminal 
          over $V=j^{-1}(U)$. The Connectivity Lemma also implies
          that ${\cO_Y}|_V = f_*\cO_E(\ulcorner -B_E \urcorner)|_V$.

\end{enumerate}
\end{proof}

\begin{lem} Let $j:Y \to (X,B)$ be an adjunction setting such that
$j(Y)$ is an exceptional component of $LCS(X,B)$. Then
the following hold:
\begin{enumerate}
    \item Let $D$ be an $\R$-Cartier divisor on $X$ not containing $j(Y)$
         in its support. Then $j:Y \to (X,B+D)$ is an adjunction setting 
         and $(B+D)_Y=B_Y+j^*(D)$.
    \item Let $\tau$ be an extraction and $j, j'$ two adjunction settings 
          making commutative the following diagram:
          \[ \xymatrix{
               Y' \ar[rr]^{\tau} \ar[dr]_{j'} &  & Y \ar[dl]^j \\
                   &                    (X,B)  &
          } \]
          Then $B_Y=\tau_*(B_{Y'})$.
    \item Let $\nu:W^\nu \to (X,B)$ be the induced adjunction setting. 
          Then $B_{W^\nu}$ is an effective divisor if $B$ is.
\end{enumerate}
\end{lem}

\begin{proof}(cf. \cite{Ka2,Ka3}) The first two properties are formal 
       consequences of the definition of the residue and 
       Remark~\ref{basic}$.6$.     
       As for the last statement, let $\nu:Y=W^\nu \to (X,B)$ be the 
       induced adjunction setting. 
       By the Kawamata-Viehweg vanishing theorem 
       $R^1\mu_*{\mathcal O}_{\tilde{X}}(-E+\ulcorner -A \urcorner)=0$, 
       hence the morphism
       $ \mu_*{\mathcal O}_{\tilde{X}}(\ulcorner -A \urcorner) 
          \to \mu_*{\mathcal O}_E(\ulcorner -B_E \urcorner) 
       $ 
is surjective. Let $Q$ be a prime divisor on $Y$ and assume by
contradiction that $b_Q<0$. Lemma~\ref{effectivity}
shows that $B_E= A|_E$ has coefficients less than $1$ over
$\eta_Q$, hence the same holds for $A$. Since $B$ is effective,
we deduce that 
$\mu_*{\mathcal O}_{\tilde{X}}(\ulcorner -A \urcorner)_Q 
=\cO_{X,Q}$. Therefore the induced map
$$
\cO_{Y,Q} \to 
 f_*{\mathcal O}_E(\ulcorner -B_E \urcorner)_Q 
$$
is surjective. This contradicts the conclusion of Lemma~\ref{effectivity}.
Our assumption was false, hence $b_Q \ge 0$.
\end{proof}

\begin{lem} Assume the Base Change Conjecture~\ref{base_change_conj} 
holds true. Let
$j:Y \to (X,B)$ be an adjunction setting such that $j(Y)$ is
an exceptional irreducible component of $LCS(X,B)$. 
Let $\sigma:Y \to W^\nu$ be the induced extraction. 
\[ \xymatrix{
       (Y,B_Y) \ar[rr]^{\sigma} \ar[dr]_j &  
                      & (W^\nu,B_{W^\nu}) \ar[dl]^\nu \\
         &                    (X,B)  &
} \]
Then
$K_Y+B_Y$ and $K_{W^\nu}+B_{W^\nu}$ are $\Q$-Cartier log divisors
and $$\sigma^*(K_{W^\nu}+B_{W^\nu})=K_Y+B_Y$$
\end{lem}

\begin{proof}
It is a formal consequence of the Base Change Conjecture applied
to $f:(E,B_E) \to W^\nu$ under the birational base change $\sigma$.
\end{proof}

\begin{rem}~\label{semin} All the concepts in this chapter are well 
defined for non-exceptional lc centers too, provided the adjunction 
calculus on  seminormal varieties is developed (representatives on 
log resolutions for lc places over $W$ are no longer irreducible, but 
they are simple normal crossings, hence seminormal). 
\end{rem}

%%%%%%%%%%%%%%%%%%%%%%%%%%%%%%%%%%%%%%%%%%
%%%%%%%%%%%%%%%%%%%%%%%%%%%%%%%%%%%%%%%%%%
\medskip
\subsection{Positivity of the moduli part}

%%%%%%%%%%%%%%%%%%%%%%%%%%%%%%%%%%%%%%%%%%
%%%%%%%%%%%%%%%%%%%%%%%%%%%%%%%%%%%%%%%%%%

Let $j:Y \to (X,B)$ be an adjunction setting, $E$ the lc place over
$j(Y)$, $f:(E,B_E) \to Y$ the induced morphism and $B_Y$ the induced 
residue on $Y$. Assume $K_Y+B_Y$ is $\R$-Cartier. This assumption 
is satisfied if $Y$ is  $\Q$-factorial and it should always hold
according to the Base Change Conjecture.
\par
Since $K_E+B_E \sim_{K,f} 0$ and $f$ is a contraction, there exists
a unique class $M_Y \in \Pic(Y) \otimes \R$ such that
$$
K_E+B_E-f^*(K_Y+B_Y) \sim_\R f^*M_Y
$$
The class $M_Y$ does not depend on the choice of the realization of $E$, and
it is called the {\em moduli part} of $K_X+B$ on $Y$ \cite{Ka2,Ka3} . An
equivalent definition of $M_Y$ is given by the following 
{\em adjunction formula}:
$$
j^*(K_X+B) \sim_\R K_Y+B_Y+M_Y 
$$

The appearance of the moduli part was first pointed out by Y. Kawamata.
It is trivial in the case $\codim(j(Y),X)=1$, according to the following 
result.
 
\begin{thm}\cite{Sh_flips} Let $j:Y \to (X,B)$ be an adjunction setting 
such that $j(Y)=W$ is a prime divisor on $X$.
Let $B_Y$ be the induced residue on $Y$.
\begin{itemize}
   \item[-] $K_Y+B_Y$ is $K$-Cartier and $r(K_Y+B_Y)$ is Cartier
            if $r(K_X+B)$ is Cartier for some $r \in \N$.
   \item[-] If $\sigma:Y \to W^\nu$ is the induced 
            extraction, then $K_Y+B_Y=\sigma^*(K_{W^\nu}+B_{W^\nu})$.
       \[ \xymatrix{
                  (Y,B_Y) \ar[rr]^{\sigma} \ar[dr]_j &  
                      & (W^\nu,B_{W^\nu}) \ar[dl]^\nu \\
         &                    (X,B)  &
       } \]
   \item[-] The moduli part $M_Y$ is trivial, that is
           $j^*(K_X+B) \sim_K K_Y+B_Y$.

   \item[-] $a(Z;X,B) \le a(j^{-1}(Z);Y,B_Y)$ 
              for every proper subvariety $Z \subset W$.
\end{itemize}
\end{thm}

\begin{proof}
Since $codim(W,X)=1$, $W$ is an exceptional lc center and
the induced morphism $f:E \to Y$ is an extraction. By 
Example~\ref{exmp_pf}$.1$, $K_Y+B_Y$ is $K$-Cartier and 
$f:(E,B_E) \to (Y,B_Y)$ is a crepant extraction. 
Therefore $\sigma:(Y,B_Y) \to (W^\nu,B_{W^\nu})$ is also crepant.
The rest is formal.
\end{proof}

\begin{thm} \label{adjthm} (cf. \cite[Theorem 1]{Ka3}) Let 
$j:Y \to (X,B)$ be an adjunction setting such that
$W=j(Y)$ is an exceptional irreducible component of $LCS(X,B)$.
Then $j$ is dominated by an adjunction setting $j':Y' \to (X,B)$ 
satisfying the following properties:
\begin{itemize}
  \item[1.] ${j'}^*(K_X+B) \sim_\Q K_{Y'}+B_{Y'}+M_{Y'}$, 
            where $B_{Y'}$ is the residue
            of $K_X+B$ on $Y'$ and $M_{Y'} \in Pic(Y') \otimes \Q$ is
            the moduli part;
  \item[2.] $M_{Y'}$ is nef;
  \item[3.] For any $\Q$-Cartier divisor on $X$ whose support does 
            not contain $W$ the following hold:
    \begin{itemize}
        \item[a)] If $(X,B+D)$ is log canonical on an open subset
                   $U \subseteq X$, then $({Y'},(B+D)_{Y'})$
                   is log canonical on ${j'}^{-1}(U)$.
        \item[b)] If $W$ is the minimal lc center of $(X,B+D)$ at a closed 
                   point $j(y) \in X$, then $(Y',(B+D)_{Y'})$ is
                   Kawamata log terminal on a neighborhood of 
                   ${j'}^{-1}(j'(y))$.
    \end{itemize}
\end{itemize} 
\end{thm}

\begin{proof}
Denote by $E$ the unique lc place over $W$. There exists an extraction
$\sigma:Y' \to Y$ such that $E \to Y$ factors through
$f:E \to Y'$ and $f:(E,B_E) \to Y'$ satisfies the property $(\star)$. 
Therefore the first two properties hold.
\par
From Lemma~\ref{stan}$.4$ and Remark~\ref{basic}$v)$, the last part 
holds for $D=0$. 
Finally, let $D$ be a $\Q$-Cartier divisor on $X$ whose
support does not contain $W$. There exists an extraction
$\tau:Y'' \to Y'$ such that 
both $f:(E,B_E) \to Y'$ and $f':(E',(B+D)_{E'}) \to Y''$ 
satisfy the property $(\star)$, and moreover $f'$ is induced
by the birational morphism $\tau:Y' \to Y$. 
By the previous step, $(Y'',(B+D)_{Y''})$ has the required properties. 
Moreover, Proposition~\ref{sigma_prop} gives the inequality
$$\tau^*(K_{Y'}+B_{Y'}+{j'}^*(D)) \le 
 K_{Y''}+B_{Y''}+(\tau \circ j')^*(D)= K_{Y''}+(B+D)_{Y''}
$$
which implies that these properties are inherited by $(Y',(B+D)_{Y'})$.
\end{proof}

For some applications, the following result of Y. Kawamata suffices.

\begin{thm}\cite[Theorem 1]{Ka3} Let $(X,B)$ be a projective log
variety and let $W$ be a minimal center of log canonical 
singularities for $(X, B)$. Assume 
there exists an effective $\Q$-divisor $B^o$ on $X$ such
that $B^o \le B$ and $(X,B)$ is Kawamata log terminal in
a neighborhood of $W$. 
\newline
Let $H$ be an ample Cartier divisor on $X$, 
and $\epsilon$ a positive rational number.
Then there exists an effective $\Q$-divisor $D_W$ on $W$
such that 
$$
(K_X + B + \epsilon H) \vert_W \sim_\Q K_W + D_W
$$
and that the pair $(W, D_W)$ is Kawamata log terminal.
\end{thm}

\begin{proof}\cite[Theorem 1]{Ka3} We assume $dimW>0$, otherwise there 
is nothing to prove. Then $a(\eta_W;X,B)=0$ and $(X,B)$ is log canonical 
in a neighborhood of $W$ according to $1.4.1$. 
Moreover, $W$ is normal from the Conectedness Lemma.
\newline
There exists an effective $\Q$-divisor $D'$ passing through $W$ such that
$W$ is an exceptional minimal lc center of $(X,B^o+(1-t)(B-B^o)+tD')$ 
for $0<t \ll 1$. 
From the previous theorem, there exists a resolution $\sigma:Y \to W$
such that the differents $(B^o+(1-t)(B-B^o)+tD')_Y$ are supported
in a simple normal divisor $Q$ for every $0<t \ll 1$ and
$$
M_t=\sigma^*(K_X+B^o+(1-t)(B-B^o)+tD'|_W)-(K_Y+B^o+(1-t)(B-B^o)+tD')_Y)
$$
is a nef divisor for every $0<t \ll 1$. 
Let $B_Y=\lim_{t \to 0} B^o+(1-t)(B-B^o)+tD'$, which is a divisor
supported in $Q$ (this is the different in the non-exceptional case). 
\newline
Then $M_0=\sigma^*(K_X+B|_W)-(K_Y+B_Y)$ is a nef divisor. Moreover,
$B_Y$ has coefficients less than $1$ and is relative effective
over $W$. Let $Q'$ be an effective $\sigma$-exceptional $\Q$-divisor 
with very small coefficients such that $-Q'+M_0+\sigma^*(\epsilon H)$ is 
ample on $Y$, hence there exists an effective $\Q$-divisor $M'$ on $Y$ such 
that $M' \sim_\Q -Q'+M_0+\sigma^*(\epsilon H)$ and $\Supp(M') \cup \Supp(Q)
\cup \Supp(B_Y)$ is a simple normal crossing divisor.
Since $K_Y+B_Y+(Q'+M') \sim_{\Q, \sigma} 0$, we obtain that
$D_W=\sigma_*(B_Y+Q'+M')$ is an effective divisor such that
$K_W+D_W$ is $\Q$-Cartier, $\sigma^*(K_W+D_W)=K_Y+(B_Y+Q'+M')$
and $K_W+D_W \sim_\Q (K_X+B)|_W$.
Finally, $(W, D_W)$ has Kawamata log terminal singularities
since it has $(Y,B_Y+Q'+M')$ as a log resolution.
\end{proof}

\begin{rem}
The same proof gives a localized version of the above theorem:
if $W$ is the minimal lc center at a closed point $x \in X$ and 
$dimW>0$, then we can choose such a divisor $D_W$ such that 
$(W,D_W)$ is Kawamata log terminal in a neighborhood of $x$.
\end{rem}

\newpage

%%%%%%%%%%%%%%%%%%%%%%%%%%%%%%%%%%%%%%%%%%%%%%%%%
%%%%%%%%%%%%%%%%%%%%%%%%%%%%%%%%%%%%%%%%%%%%%%%%%
%%%%%%%%%%%%%%%%%%%%%%%%%%%%%%%%%%%%%%%%%%%%%%%%%

\section{Building singularities}

%%%%%%%%%%%%%%%%%%%%%%%%%%%%%%%%%%%%%%%%%%%%%%%%%
%%%%%%%%%%%%%%%%%%%%%%%%%%%%%%%%%%%%%%%%%%%%%%%%%
%%%%%%%%%%%%%%%%%%%%%%%%%%%%%%%%%%%%%%%%%%%%%%%%%

Let $(X,B)$ be a log pair, $H$ an ample $\Q$-Cartier divisor and 
$x \in X \setminus LCS(X,B)$ a closed point. 
For every $c>0$ denote by $\cS_x(B,H;c)$ the set of all effective 
$\Q$-Cartier divisors $D$ on $X$ such that
\begin{itemize}
   \item[-] $D \sim_\Q c H$;
   \item[-] $a(x;B+D)=0$.
\end{itemize}
Since $H$ is ample, $\cS_x(B,H;c) \ne \emptyset$ for $c$ sufficiently
large, so the following infimum is well defined: 
$$
bld_x(B;H)=\inf\{c>0; \cS_x(B,H;c) \ne  \emptyset \}
$$
The problem of effective building of singularities consists
of finding upper bounds for $bld_x(B;H)$ in terms of
the invariants of the local singularity $x \in (X,B)$ 
and the global numeric properties of $H$.

\begin{example} 
   \begin{itemize}
      \item[a)] If $X$ is a curve and $H$ is a $\Q$-ample divisor then
$$bld_x(B;H)=a(x;X,B)/deg_X(H)$$
      \item[b)] Let $x$ be the vertex of the singular conic  
 $X \subset {\mathbb P}^3$ and let $H$ be the hyperplane section. Then 
 $bld_x(0;H)=a(x;X,0)=1$. Note that $H^2=2$ and 
 $\inf_{C \subset X} (H \cdot C)=1$.
  \end{itemize}
\end{example}

\begin{defn}
\begin{itemize}
   \item[a)] For a nef divisor $H \in Div(X) \otimes \R$ on a complete
variety $X$ we denote by 
$deg_X(H)=(H^{dimX})_X \in \R$ 
its top self intersection.
   
   \item[b)] We say that $H$ is {\em normalized at a closed point} 
             $x \in X$ if $deg_W(H|_W) \ge 1$ for every closed subvariety 
             $x \in W \subseteq X$.

   \item[c)] (V. V. Shokurov) A divisor $H$ on a variety 
$X$ {\em has height at least $h$} if $$H \equiv \sum_j h_j H_j +N$$ 
where $H_j$ are ample Cartier divisors, $h_j >0 \ \forall j$, 
$\sum_j h_j \ge h$, and $N$ is a nef $\R$-divisor.
\par
Note that $\frac{1}{h} H$ is normalized at every closed point $x \in X$
if $H$ has height at least $h$.
 \end{itemize}
\end{defn}

Inspired by the curve case we expect the following conjecture
to hold:

\begin{conj}\label{bld_conj} Let $(X,B)$ be a log variety,
$x \in X\setminus LCS(X,B)$ a closed point and $H \in Div(X) \otimes \Q$ 
an ample divisor normalized at $x$. Then
$$
bld_x(B;H) \le a(x;X,B)
$$
\end{conj}

%%%%%%%%%%%%%%%%%%%%%%%%%%%%%%%%%%%%%%%%%%%%%%%%%
%%%%%%%%%%%%%%%%%%%%%%%%%%%%%%%%%%%%%%%%%%%%%%%%%
\medskip
\subsection{A first estimate}

%%%%%%%%%%%%%%%%%%%%%%%%%%%%%%%%%%%%%%%%%%%%%%%%%
%%%%%%%%%%%%%%%%%%%%%%%%%%%%%%%%%%%%%%%%%%%%%%%%%

We use the standard abbreviation $h^0(X,\cN)=\dim_k H^0(X,\cN)$
for the dimension of the space of global sections of a coherent 
sheaf $\cN$ on a variety $X$.

\begin{defn}\label{generic_contr} We say that a morphism of varieties 
$f:X \to P$ is
a {\it generic contraction} if it is proper and the induced 
map $\cO_{f(X)} \to f_*\cO_X$ is an isomorphism in the generic point
of the induced subscheme $f(X) \subset P$. 
\end{defn}
The point of the Definition~\ref{generic_contr}, formalized
in Lemma~\ref{almostall}, is that if $f:X \to P$
is a generic contraction and $\cH$ is an ample line bundle
on $P$, even if not all global sections of $f^*\cH$ are pull backs of
global sections of $\cH$, almost all of them are. 
For our purposes, $X$ will be the normalization or a resolution of a 
generically reduced subvariety of $P$.

\begin{lem}\label{almostall}
Let $f:X \to P$ be a generic contraction, let $\cH$ be an ample line 
bundle on $P$, and let
$$V_k=\im[H^0(P,\cH^{\otimes k}) \to H^0(X,f^*\cH^ {\otimes k})]$$ 
Then
$$
\lim_{k \to \infty} \frac{h^0(X,f^*\cH^{\otimes k})}{\dim V_k}=1.
$$
\end{lem}

\begin{proof} Let $\cG=\Coker(\cO_P \to f_*\cO_X)$.
The cohomological interpretation of ampleness gives the exactness of the
following sequence
$$
H^0(P,\cH^{\otimes k}) \to H^0(P,f_*\cO_X \otimes \cH^{\otimes k})\to
 H^0(P,\cG \otimes \cH^{\otimes k}) \to 0
$$
for $k \in \N$ large enough. 
Note that $H^0(X,f^*\cH^{\otimes k})=H^0(P,f_*\cO_X \otimes \cH^{\otimes k})$ 
by the projection formula.
There exist 
polynomials $P(k)$ and $Q(k)$, of degrees $\dim \Supp(f_*\cO_X)$ and 
$\dim \Supp(\cG)$ respectively, such that
$$h^0(P,f_*\cO_X \otimes \cH^{\otimes k})=
\chi(P,f_*\cO_X \otimes \cH^{\otimes k})=P(k),$$ 
$$h^0(P,\cG \otimes \cH^{\otimes k})=
\chi(P,\cG \otimes \cH^{\otimes k})=Q(k)$$
for $k \in \N$ large enough.
Note that $\dim \Supp(f_*\cO_X)>\dim \Supp(\cG)$ since $f$ is a generic 
contraction. Therefore
$
\dim V_k=h^0(P,f_*\cO_X \otimes \cH^{\otimes k})-
h^0(P,\cG \otimes \cH^{\otimes k})
$
is a polynomial in $k$ for $k$ large enough. 
Moreover, it has the same degree and leading coefficient as $P(k)$,
hence the claim.
\end{proof}

\begin{lem}\label{mult} \cite[1.3]{Sh_nv}
Let $f:X \to P$ be a generic contraction of normal varieties and
$H$ an ample $\Q$-Cartier divisor on $P$ such that
$deg_X(f^*H)>1$.
There exists a natural number $k \in \N$ such that 
$kH$ is an integer divisor and for any nonsingular point 
$x \in X$ there exists an effective divisor $D_x \in |kH|$ whose
support does not contain $f(X)$ and $\mult_x(f^*D_x)>k$. 
\end{lem}

\begin{proof} We use the notation from the previous lemma.
Denote also $d=dimX$ and $q^d=deg_X(f^*H)$.
We have to prove that the natural map of vector spaces
$$
\varphi_k: V_k \to H^0(X,f^*\cO(kH) \otimes \cO_X/m^{kq}_x)
$$
has a nontrivial kernel for $k \in \N$ sufficiently large and divisible. 
Since $f^*H$ is nef and big on $X$, we have 
$$\dim H^0(X,f^*\cO(kH))=\frac{q^d}{d!} k^d+O(k^{d-1})$$
for $k \in \N$ sufficiently large and divisible.
From the previous lemma, we deduce that
$$
\dim V_k > \frac{q^d}{d!} k^d+O(k^{d-1})
$$
for $k \in \N$ large enough. 
However, 
$$\dim H^0(X,f^*\cO(kH) \otimes \cO_X/m^{kq}_x)
=\dim H^0(X,\cO_X/m^{kq}_x)=\frac{q^d}{d!} k^d + O(k^{d-1})$$
 hence
the morphism $\varphi_k$ cannot be injective for
$k$ large enough and divisible. Notice that
this $k$ is independent of the choice of the smooth point $x$.
\end{proof}

%%%%%%%%%%%%%%%%%%%%%%%%%%%%%%%%
%%%%%%%%%%%%%%%%%%%%%%%%%%%%%%%%

\medskip
\subsection{A quadratic bound}

%%%%%%%%%%%%%%%%%%%%%%%%%%%%%%%%
%%%%%%%%%%%%%%%%%%%%%%%%%%%%%%%%

\begin{thm}(cf. \cite[Theorem 6.7.1]{Ko2}) ~\label{Kollar}
Let $f:X \to P$ be a generic contraction of normal 
varieties, $x \in X$ a closed point, and let $P$ be polarized by 
a $\Q$-ample divisor $H$ such that $deg_X(f^*H)>1$. 
Then there exists an effective divisor 
$D_x \in Div(P) \otimes \Q$ satisfying
the following properties:
\begin{itemize}
  \item[1.] $D_x \sim_\Q c H$ for some rational number $c< dimX$ and
            $f(X)$ is not contained in the support of $D_x$,
  \item[2.] $x \in LCS(X,B+f^*D_x)$ for every effective $\R$-Weil divisor
            $B$ on $X$ such that $K_X+B$ is $\R$-Cartier.
\end{itemize}
\end{thm}

\begin{proof} The above statement is stronger than 
\cite[Theorem 6.7.1]{Ko2}, 
but with the same proof, presented here for completeness.
\newline
{\it Step 1.} Assume first that $x \in X$ is
nonsingular. Lemma~\ref{mult} gives a divisor
$D'_x \sim kH$ such that $B'$ does not contain $f(X)$ in its support 
and 
$mult_x(f^*D'_x) >k$. Then $B_x=\frac{dimX}{k}D'_x$ satisfies the 
required 
properties at $x$. Indeed, if $K_X+B$ is $\R$-Cartier, $B$ is effective 
and $E$ is the exceptional divisor on the blow-up of $X$ in $x$, then
$$
a(E;B+B_x)=dimX-\mult_x(B)-\frac{\mult_x(f^*D'_x)}{k}dimX <0.
$$
{\it Step 2.} We are left with the case when $x \in X$ is a singular
point. Since the integer $k$ does not depend on the choice
of the smooth point, we can assume there exists a smooth pointed
curve $(C,0)$, a morphism $g:C \to X$ such that $g(0)=x$, 
$g(C \setminus \{0\}) \subset Reg(X)$ and a morphism 
$\tilde{g}:C \to |kH|$
such that $B_{g(c)}=\frac{dimX}{k}\tilde{g}(c)$ satisfies the
required property at $g(c) \in X$ for each $c \in C \setminus \{0\}$.
By Proposition~\ref{variation}, $B_x=D_{g(0)}$ also satisfies the
required property at $g(0)=x$.
\end{proof}

The following lemma extends Theorem~\ref{Kollar} to relative
effective pseudo-boundaries.

\begin{lem}\label{releff} Let $\sigma:Y \to W$ be a birational 
contraction, $w \in W$ a closed point, $K_Y+B_Y$ a log divisor and 
$D$ an 
effective $\R$-Cartier on $W$ with the following properties:
\begin{itemize}
  \item[1.] $K_Y+B_Y+M \sim_{\Q,\sigma} 0$, where 
                $M \in Div(Y) \otimes \Q$ is a nef divisor,
  \item[2.] $\sigma_*B_Y$ is an effective divisor,
  \item[3.] $w \in LCS(W,B+D)$ for every effective $\R$-divisor 
            $B \ge \sigma_*B_V$ such that $K_W+B$ is $\R$-Cartier.
\end{itemize}
Then $\sigma^{-1}(w) \cap LCS(Y,B_Y+\sigma^*D) \ne \emptyset$.
\end{lem}

\begin{proof}(cf. \cite{Ka3}) The assumptions are invariant under 
blow-ups on $Y$, so we can assume that 
$S=Exc(\sigma) \cup \Supp(B_Y) \cup \Supp(\sigma^*D)$ has normal 
crossing.
By contradiction, there exists $0 < \epsilon \ll 1$ such that
$$a(E;B_Y+\sigma^*D) \ge \epsilon $$ for every prime divisor $E$ on $Y$
such that $\sigma^{-1}(w) \cap E \ne \emptyset$.
Let $A$ be an effective exceptional divisor on $Y$ with coefficients
less than $\epsilon$, such that $-A$ is $\sigma$-ample and
$M-A+\sigma^*M$ is $\Q$-ample for some ample divisor $H$ on $W$.
Let $M' \sim_\Q M-A+\sigma^*M$ be an effective divisor with 
coefficients
less than $\epsilon$ such that $S \cup Supp(M')$ has normal 
crossing.
In particular, $a(E;B_Y+\sigma^*D+A+M')=
a(E;B_Y+\sigma^*D)-\mult(E;A+M')>0$
for every prime divisor $E$ on $Y$ such that 
$\sigma^{-1}(w) \cap E \ne \emptyset$. Therefore
$$
\sigma^{-1}(w) \cap LCS(Y,B_Y+\sigma^*D+A+M')= \emptyset
$$
Let $B_W= \sigma_*(B_Y+A+M')$. Since $K_Y+B_Y+A+M' \sim_{\Q,\sigma} 0$, 
we deduce that $K_W+B_W$ is $\R$-Cartier and 
$$K_Y+B_Y+A+M'=\sigma^*(K_W+B_W)$$  
Moreover, $B_W \ge \sigma_*B_Y$, hence 
$\sigma^{-1}(w) \cap LCS(Y,B_Y+\sigma^*D+A+M')\ne \emptyset$.
Contradiction.
\end{proof}

\begin{prop}\label{indstep} Let $(X,B)$ be log variety and 
$x \in X$ a closed point such that $a(x;X,B) \ge 0$ and $(X,B)$ 
has a normalized minimal lc center $W$ at $x$.
Let $H$ be a polarization of $X$ such that $deg_W(H|_W)>1$. 
Assume there exists an effective
divisor $B^o<B$ such that $(X,B^o)$ is Kawamata log terminal at $x$.
\par
Then there exists an effective divisor $B_1 \in Div(X) \otimes \Q$
such that
\begin{itemize}
   \item[a)] $B_1 \sim_\Q c_1 H$, $c_1< dimW_1$
   \item[b)] $a(x;B+B_1) \ge 0$ and $(X,B+B_1)$ has a normalized
             minimal lc center $W_1$ at $x$. 
   \item[c)] $W_1 \subset W$ and $dimW_1<dimW$
\end{itemize}
\end{prop}

\begin{proof}
Let $Y \to W$ be a resolution as in Theorem~\ref{adjthm}, and
let $\sigma^\nu:Y \to W^\nu$ be the induced map to the normalization
of $W$. Since $W$ is normal at $x$, we can harmlessly say $x \in W^\nu$.
Let $B_x \sim_\Q c H\ (c < dimW)$ be the divisor obtained by 
applying Theorem~\ref{Kollar} to $x \in W^\nu \to (X,H)$. From 
Lemma~\ref{releff} applied to $Y \to W^\nu$ and $D=B_x$, we deduce that 
$LCS(Y,B_Y+\sigma^*(B_x|_W)) \cap \sigma^{-1}(x) \ne \emptyset$, 
that is $K_Y+B_Y+\sigma^*(B_x|_W)$ is not klt near $\sigma^{-1}(x)$.
\par
Let $\alpha=\sup\{t>0; (X,B+tB_x) \textnormal{ is log canonical at } x\}$.
Theorem ~\ref{adjthm}$3.b)$ implies that $\alpha \le 1$. 
It is clear that $(X,B+\alpha B_x)$ is log canonical at $x$, with minimal 
lc center $W_1$ at $x$, strictly included in $W$. After a small 
perturbation of $B_x$ along $H$, we can assume that $W_1$ is normalized
too. Therefore $B_1=\alpha B_x$ and $c_1=\alpha c$ have the 
required properties.
\end{proof}

\begin{thm}\label{mainbld} \cite[Theorem 6.4]{Ko2}
Let $(X,B)$ be a log variety, $x \in X\setminus LCS(X,B)$ a closed point 
and $H \in Div(X) \otimes \Q$ an ample divisor normalized at $x$.
Then $$bld_x(B;H) < \frac{1}{2}dimX(dimX+1)$$
\end{thm}

\begin{proof}
Set $W_0=X$. By Proposition ~\ref{indstep},
there exists an effective $\Q$-Cartier divisor $B_1$ such that 
\begin{enumerate}
   \item $B_1 \sim_\Q c_1 H$, $c_1+\epsilon < dimW_0$ for
             some small enough $\epsilon$
   \item $(X,B+B_1)$ is log canonical at $x$ with normalized
             and minimal lc center $W_1$ at $x$
   \item $dimW_1<dimW_0$ 
\end{enumerate}
We repeat the previous step for $W_1$ and so forth, only that we
apply Proposition~\ref{indstep} for 
$(1+\frac{\epsilon}{m})H$ instead of $H$, where 
$m=\frac{dimX(dimX+1)}{2}$. 
Thus we obtain a sequence of divisors $B_1,B_2, \cdots$ 
such that
\begin{enumerate}
   \item $B_j \sim_\Q c_j H$
   \item $c_1+\epsilon < dimW_0$ and 
$c_j<(1+\frac{\epsilon}{m})dimW_{j-1}$ 
         for every $j \ge 2$
   \item $(X,B+\sum_{k=1}^j B_k)$ is log canonical at $x$ 
with normalized
             and exceptional lc center $W_k$ at $x$
   \item $dimW_{j+1}<dimW_j$ for every $j \ge 0$ 
\end{enumerate}
We stop this recursive process at some step $s$ for which 
$W_s=\{x\}$.
This definitely happens for some $s \le dimW_0$, due to property 
$4$ above.
We have $\sum_{k=1}^s c_k< dimW_0-\epsilon +
   (1+\frac{\epsilon}{m})\sum_{k=2}^sdimW_{k-1} \le 
\sum_{k=0}^s dimW_c \le
m$. Therefore $B_x=\sum_{k=1}^s B_k \sim_\Q cH$, $c<m$ and 
$a(x;B+B_x)=0$.
\end{proof}

%%%%%%%%%%%%%%%%%%%%%%%%%%%%%%%%%%%%%%%%%%%%%%%%%
%%%%%%%%%%%%%%%%%%%%%%%%%%%%%%%%%%%%%%%%%%%%%%%%%
\medskip
\subsection{The conjectured optimal bound}

%%%%%%%%%%%%%%%%%%%%%%%%%%%%%%%%%%%%%%%%%%%%%%%%%
%%%%%%%%%%%%%%%%%%%%%%%%%%%%%%%%%%%%%%%%%%%%%%%%%

We discuss in the section the connection between
Conjecture~\ref{bld_conj} and the Adjunction Conjecture.
 
\begin{conj}\label{bld_step} Let $(X,B)$ be a log pair which is 
Kawamata log terminal in a neighborhood of a closed point $x \in X$.
Let $H \in Div(X) \otimes \Q$ be an ample divisor 
normalized at $x$ and fix $0<\epsilon \ll 1$. 
Then there exists an effective $\Q$-Cartier divisor $D$ with the 
following properties:
\begin{itemize}
  \item[1.] $D \sim_\Q c H$
  \item[2.] $0<(1-\epsilon) c \le a(x;B)-a(x;B+D)$
  \item[3.] $(X,B+D)$ is maximally log canonical at $x$
\end{itemize}
\end{conj}

\begin{lem} Assume the first two properties of the different
in the Adjunction Conjecture hold true. 
Then the two conjectures ~\ref{bld_step} and ~\ref{bld_conj} 
are equivalent.
\end{lem}

\begin{proof} We just need to show that Conjecture~\ref{bld_step}
implies Conjecture ~\ref{bld_conj}. We use induction on $dimX$.
Fix $0<\epsilon \ll 1$ and let $D_0 \sim_\Q c_0 H$ be a divisor given by 
Conjecture~\ref{bld_step}. We may assume that $W$, the minimal lc center 
of $(X,B+D_0)$ at $x$, is normalized. In particular, $x$ is a normal point 
of $W$.
\newline
If $W=\{x\}$, then $(1-\epsilon)c_0 \le a(x;X,B)-a(x;X,B+D_0)=a(x;X,B)$. 
We hence assume that $dimW>0$.
Let $(B+D_0)_{W^\nu}$ be the different of $K_X+B+D_0$ on $W^\nu$. 
Then $(W^\nu,(B+D_0)_{W^\nu})$ is klt at $x$ and 
$H_{W^\nu}=H|_{W^\nu}$ is ample normalized at $x$. By induction, there
exists an effective divisor $D_1 \in Div(X) \otimes \Q$ with the following
properties:
 \begin{enumerate}
    \item $D_1 \sim_\Q c_1 H$
    \item $(1-\epsilon)c_1 \le a(x;W^\nu,(B+D_0)_{W^\nu})$ 
    \item $a(x;W^\nu,(B+D_0)_{W^\nu}+{D_1}|_{W^\nu})=0$. 
  \end{enumerate}
By precise inverse of adjunction, we deduce that $a(x;X,B+D_0+D_1)=0$.
Therefore $D=D_0+D_1 \in \cS_x(B,H;c)$, with $c=c_0+c_1$. But
$(1-\epsilon)c=(1-\epsilon)c_0+(1-\epsilon)c_1\le
   a(x;B)-a(x;B+D_0)+a(x;W^\nu,(B+D_0)_{W^\nu})$.
By precise inverse of adjunction again, $a(x;W^\nu,(B+D_0)_{W^\nu}) \le 
    a(x;B+D_0)$, hence 
$$
(1-\epsilon)c \le a(x;B).
$$ 
Letting $\epsilon$ tend to $0$ we deduce that $bld_x(B;H) \le a(x;B)$
\end{proof}

\begin{lem} Conjecture~\ref{bld_step} holds true if either $x \in X$ 
is a nonsingular point, or $H$ is a base point free ample Cartier 
divisor. 
\end{lem}
\begin{proof}
We first show that there exists an effective $\Q$-Cartier divisor 
$D_1 \sim_\Q H$ such that
$$
\inf_F  
      \mult_F(f^*D_1)  \ge 1-\epsilon
$$
where the infimum is taken after all prime divisors $F$ on 
birational extractions $f:Y \to X$ with $f(F)=\{x\}$.
Indeed, since $deg_X(H) > (1-\epsilon)^{dimX}$, Lemma~\ref{mult} gives 
the required divisor if $x$ is a nonsingular point. If $H$ is a base 
point free ample Cartier divisor, then we may take $D_1$ to be any 
divisor passing through $x$ in the linear system $|H|$.
\par
 Let $c>0$ such that $(X,B+cD_1)$ is maximally log canonical at $x$. 
We can assume that the minimal lc center at $x$ is normalized. 
\newline
We claim that $D=cD_1$ has the desired properties.
Let $f:Y \to X$ be an extraction and $F \subset Y$ a prime divisor such 
that $f(F)=\{x\}$. Then
$a(F;B+D)=a(F;B)-c \cdot \mult_F(f^*D_1) \le 
a(F;B)-(1-\epsilon)c$.
Therefore
$$
a(F;B+D) \le a(F;B)-(1-\epsilon)c
$$
for every prime divisor $F$ with center $\{x\}$ on $X$. Taking infimum
after all such $F's$ we get the desired inequality.
\end{proof}

In particular, Conjecture~\ref{bld_conj} follows from the precise inverse
of adjunction if $H$ is a very ample Cartier divisor.

%%%%%%%%%%%%%%%%%%%%%%%%%%%%%%%%%%%%%%
%%%%%%%%%%%%%%%%%%%%%%%%%%%%%%%%%%%%%%
\medskip
\subsection{Global generation of adjoint line bundles}

%%%%%%%%%%%%%%%%%%%%%%%%%%%%%%%%%%%%%%
%%%%%%%%%%%%%%%%%%%%%%%%%%%%%%%%%%%%%%

The main application of effective building of isolated log canonical 
singularities is to the global generation of adjoint line bundles.
The key step is the following result of Y. Kawamata.

\begin{prop}\label{keystep}\cite[Proposition 2.3]{Ka1} 
Let $(X,B)$ be a log variety, $x \in X \setminus LCS(X,B)$
a closed point and $H$ an ample $\Q$-Cartier divisor.
Assume $L$ is a Cartier divisor on $X$ such that
\begin{itemize}
  \item[-] $L \equiv K_X+B+h H$; 
  \item[-] $h>bld_x(B;H)$.
\end{itemize}
Then the sheaf $\cI(X,B) \otimes \cO_X(L)$ is generated by 
global sections at $x$. This means that there exists a global 
section 
$s \in H^0(X,\cO_X(L))$ such that $s|_{LCS(X,B)}=0$ and 
$s(x) \ne 0$.
\end{prop}

\begin{proof} Let $D \in \cS_x(B;H)$ with $D \sim_\Q cH$ and 
$c<h$.
According to Lemma~\ref{gg}, we may assume that $\{x\}$ 
is a normalized lc center, that is $LCS(X,B+D) \cap U=\{x\}$
for some neighborhood $U$ of $x$.
Since $L \equiv K_X+B+D+(h-c)H$ and $(h-c)H$ is ample, the 
extension of Kawamata-Viehweg vanishing gives the surjection
$$
H^0(X,L) \to H^0(LCS(X,B+D), L|_{LCS(X,B+D)}) \to 0
$$
But $LCS(X,B+D)=\{x\} \sqcup X'$, where $X'$ is a closed subscheme of 
$X$.
In particular, 
$$
H^0(LCS(X,B+D), L|_{LCS(X,B+D)})=
H^0(\{x\},L|_{\{x\}}) \oplus H^0(X',L|_{X'})
$$
The lifting $s$ of the global section $(1,0)$ satisfies the required 
properties. Clearly, $s(x) \ne 0$, and $s|_{LCS(X,B)}=0$ since 
$LCS(X,B)$ is a closed subscheme of $X'$.
\end{proof}

\begin{lem} For any divisor $D \in \cS_x(B,H;c)$ there exists an
effective $\Q$-Cartier divisor $D' \sim_\Q H$ 
such that $D_\epsilon=(1-\epsilon)D+\epsilon D' \in 
\cS_x(B,H;(1-\epsilon)c+\epsilon)$ 
and $$LCS(X,B+D_\epsilon)\cap U=\{x\}$$ for
some open neighborhood $U$ of $x$ and for all $0<\epsilon \ll 1$.
\end{lem}

\begin{cor}\label{gg} Let $(X,B)$ be a log variety, 
$x \in X\setminus LCS(X,B)$ a closed point, and 
$H \in Div(X) \otimes \Q$ 
an ample divisor normalized at $x$. Assume 
$L$ is a Cartier divisor on $X$ such that $L \equiv K_X+B+h H$ for
some real number $h>0$. Then
the sheaf $\cI(X,B) \otimes \cO_X(L)$ is generated by 
global sections at $x$ if one of the following holds:
\begin{itemize}
   \item[-]\cite{AS,Ko2} $h > \frac{dimX(dimX+1)}{2}$
   \item[-] $h > a(x;X,B)$ if Conjecture~\ref{bld_conj} holds true.
\end{itemize}
\end{cor}

\begin{conj}(T. Fujita) Let $L$ be an ample line bundle
on a projective algebraic variety $X$. Then the adjoint
line bundles $K_X+mL$ are generated by global sections for 
$m>dimX$.
\end{conj}

\begin{rem} Fujita's Conjecture is implied by Conjecture~\ref{bld_conj}.
Indeed, set $B=0$ and $H=L$. Note that $LCS(X,B)=\emptyset$ and
$H$ is normalized in any closed point $x \in X$, since $H$ is a 
Cartier divisor. For any closed point $x \in X$ we have $a(x;X,0)=dimX$, 
hence Conjecture~\ref{bld_conj} would imply that $L$ is globally 
generated at $x$ for $m > dimX$.
\end{rem}

\newpage

\addcontentsline{toc}{section}{\bf References}   

%%%%%%%%%%%%%%%%%%%%%%%%%%%%%%%%%%%%%%%
%%%%%%%%%%%%%%%%%%%%%%%%%%%%%%%%%%%%%%%

\vspace{1in}

DEPARTMENT OF MATHEMATICS,
THE JOHNS HOPKINS UNIVERSITY,
3400 NORTH CHARLES STREET,
BALTIMORE MD 21218\\
email: ambro$@$chow.mat.jhu.edu.

\end{document}